\documentclass{amsart}
\usepackage{amssymb,amsmath,amsthm}
\usepackage[dvips]{graphicx}

\let\oldmarginpar\marginpar
\renewcommand\marginpar[1]{\-\oldmarginpar[\raggedleft\footnotesize #1]%
{\raggedright\footnotesize #1}}

\begin{document}

\newtheorem{theorem}{Theorem}[section]
\newtheorem{corollary}[theorem]{Corollary}
\newtheorem{lemma}[theorem]{Lemma}
\newtheorem{proposition}[theorem]{Proposition}
\theoremstyle{definition}
\newtheorem{definition}[theorem]{Definition}
\theoremstyle{remark}
\newtheorem{remark}[theorem]{Remark}
\theoremstyle{definition}
\newtheorem{example}[theorem]{Example}

\def\R{{\mathbb R}}
\def\H{{\mathbb H}}
\def\rank{{\text{rank}\,}}
\def\bd{{\partial}}

\title{Projective Deformations of Hyperbolic Coxeter 3-Orbifolds}
\author{Suhyoung Choi}
\address{Department of Mathematical Sciences, KAIST, Daejeon 305-701, Republic of Korea}
\email{shchoixk@gmail.com}
\thanks{The research of the first author was supported by the National Research Foundation of Korea (NRF) grant funded by the Korea government (MEST)(No. 2009-0057445).}

\author{Craig D. Hodgson}
\address{Department of Mathematics and Statistics, University of Melbourne, Victoria 3010, Australia}
\email{craigdh@unimelb.edu.au}
\thanks{The research of the second author was partially supported by Australian Research Council grant DP0663399.}

\author{Gye-Seon Lee}
\address{Department of Mathematical Sciences, KAIST, Daejeon 305-701, Republic of Korea}
\email{smileabacus@kaist.ac.kr}
\thanks{The research of the third author was supported by Basic Science Research Program through the National Research Foundation of Korea (NRF) funded by the Ministry of Education, Science and Technology (KRF-2008-621-C00003).}






\begin{abstract}
By using Klein's model for hyperbolic geometry, hyperbolic structures on orbifolds or manifolds
provide examples of real projective structures.
By Andreev's theorem, many 3-dimensional reflection orbifolds admit a
finite volume hyperbolic structure, and such a hyperbolic structure is unique.
However, the induced real projective structure on some such 3-orbifolds deforms into
a family of real projective structures that are not induced from hyperbolic structures.
In this paper, we find new classes of compact and complete hyperbolic reflection 3-orbifolds
with such deformations. We also explain numerical and exact results on projective deformations of some compact hyperbolic cubes and dodecahedra.
\end{abstract}

\maketitle
\section{Introduction}\label{intro}
\addcontentsline{toc}{section}{Introduction}

A smooth \emph{$n$-dimensional orbifold} is a Hausdorff space locally modelled on quotients of open subsets of $\mathbb{R}^n$ by finite groups of diffeomorphisms (see \cite{Thurston,BHaf,Choi2004} for detailed discussions).
 In this paper,  we deal with \emph{good} orbifolds, which are quotients of a manifold by a discrete group acting properly discontinuously, perhaps with fixed points. The image of the fixed points sets of non-trivial group elements forms the \emph{singular locus} of the orbifold.
Isomorphisms of such orbifolds are just diffeomorphisms conjugating the discrete group actions.

Given a Lie group $G$ acting transitively and effectively on an $n$-dimensional manifold $X$,
Ehresmann introduced  the idea of a
$(G, X)$-structure on an $n$-orbifold as locally modelling the orbifold on open subsets of $X$
modulo finite subgroups of $G$, with transition maps given by elements of $G$.
We refer to \cite{Thurston,Thurston_book,Choi2004} for the details.

When an orbifold $M$ admits such a $(G,X)$-structure, Thurston \cite{Thurston} showed that there exists a simply connected manifold $\tilde{M}$ and a discrete group $\Gamma$ of deck transformations so that the quotient orbifold $\tilde{M}/\Gamma$ is isomorphic to $M$.
Then $\tilde{M}$ is said to be a {\em universal cover} of $M$ and $\Gamma$ is the {\em orbifold fundamental group} of $M$; these are determined uniquely up to diffeomorphism and isomorphism respectively.
We write $\pi_1(M)=\Gamma$.   (Note that there is also a definition of orbifold fundamental group by Haefliger using paths, see \cite{BHaf}.)

Given a $(G,X)$-structure on an $n$-orbifold $M$, we can define an immersion
$D$ from the universal cover $\tilde{M}$ to $X$
and a homomorphism
\begin{displaymath}
h\; :\; \pi_1(M) \rightarrow G,
\end{displaymath}
where $\pi_1(M)$ denotes the orbifold fundamental group of $M$.
Here $D$ is called a \emph{developing map} and $h$ a \emph{holonomy homomorphism}
for the $(G,X)$-structure on $M$, and $D$ satisfies the equivariance condition
$$
D(\gamma \cdot m) = h(\gamma) \cdot D(m), \text{ for all } \gamma \in \pi_1(M), m \in \tilde M
.$$
Note that $(D,h)$ is determined only up to the following action:
\begin{eqnarray} \label{PGLaction}
(D,h(\cdot)) \mapsto (g\circ D, g\circ h(\cdot)\circ g^{-1})
\end{eqnarray}
for $g \in G$.
Conversely, the development pair $(D,h)$ determines the $(G,X)$-structure.
(See Thurston \cite[chap. 3]{Thurston_book}.)

When $M$ is a closed orbifold,
we define $\tilde{\mathfrak{D}}(M)$ to be the space of equivalence classes of
development pairs of
$(G,X)$-structures on $\tilde{M}$  modulo isotopies of $\tilde{M}$ commuting with the deck transformation group.
Here, the space of development pairs is equipped with the $C^1$-topology and $\tilde{\mathfrak{D}}(M)$ is endowed with the quotient topology.
The \emph{deformation space} $\mathfrak{D}(M)$ of
$(G,X)$-structures on the orbifold $M$
is
the quotient space of $\tilde{\mathfrak{D}}(M)$ by the action of $G$ given in equation (\ref{PGLaction}).
(See \cite{Choi2004}, \cite{Choi2006}.)
We can also think of ${\mathfrak{D}}(M)$ as the space of $(G,X)$-structures on $M$
up to the equivalence relation given by isotopy in $M$.

When $M$ is non-compact or has boundary, $\tilde{\mathfrak{D}}(M)$ and ${\mathfrak{D}}(M)$ are defined as spaces
of $(G,X)$-structures on $M$  up to the equivalence relation given by isotopy and ``thickening'' of the geometric structure near the ends or boundary of $M$;
see \cite{Thurston,CEG}.
One of the authors of this paper is writing a more complete version of this theory in \cite{Choi2010}.

In this paper we study real projective structures and hyperbolic structures.
{\em Real projective geometry} is given by the group $PGL(n+1,\R)$ acting by projective
transformations on the projective space $\R P^n$. We can represent {\em hyperbolic geometry}
using Klein's projective model:  hyperbolic space is an open ball $B$ in $\R P^n$,
and the group of hyperbolic isometries is the subgroup $PO(1,n)$  of
$PGL(n+1,\mathbb{R})$ preserving $B$. Hence hyperbolic orbifolds and manifolds naturally
have induced real projective structures.

The study of real projective structures was originally introduced by E. Cartan, and
was continued by many people including Chern, Kuiper, Koszul, Milnor, and Benzecri in the late 1950's.
In the 1960's, it was unknown whether every real projective structure arises from a hyperbolic structure.
In $1967$, Kac and Vinberg \cite{Kac1967} discovered real projective reflection orbifolds that are not hyperbolic.
Later, the theory of real projective structures on $2$-dimensional manifolds and orbifolds
was developed in Goldman's senior thesis at Princeton from 1977 written under Thurston, and by  Goldman and Choi \cite{Goldman1990}, \cite{cdcr}, \cite{Choi2005} from the 1990's onwards.
In addition, Cooper, Long and Thistlethwaite \cite{Cooper2006} investigated
whether the closed hyperbolic 3-manifolds of the Hodgson-Weeks census could be deformed into other real projective structures.

We will focus on $3$-dimensional reflection orbifolds whose underlying space is
homeomorphic to a $3$-dimensional convex polyhedron, and whose singular locus is its boundary
(made up of mirrors).
The fundamental group of such an orbifold is a \emph{Coxeter group}, i.e.
a group with a set $\{r_1, \ldots, r_m \}$ of generators and the following set of defining relations:
{\setlength\arraycolsep{2pt}
\begin{eqnarray*}
r_i^2&=&1 \text{ for all } i,\\
(r_i r_j)^{n_{ij}}&=&1 \text{ for some $i$ and $j$ with $n_{ij}=n_{ji} \geq 2$}.
\end{eqnarray*}}
Here, $r_i$ represents a reflection in the $i$th silvered face of $P$,
and $r_i r_j$ represents a rotation of order $n_{ij}$
about an edge  where the $i$th and $j$th faces meet. The stabilizer of each face is the group
$\mathbb{Z}_2$ generated by reflection in the face, and the stabilizer of each edge
is the dihedral group $D_{n_{ij}}$ generated by reflections in the adjacent faces.

Let $P$ be a fixed $3$-dimensional convex polyhedron, and
assign an order $n_e \ge 2$ to each edge $e$ of $P$.
If any vertex of $P$ has more than three edges incident, or has orders of the incident edges not of the form
\begin{displaymath}
(2,2,k) \:\:\text{with}\: k\geq 2,\:\: (2,3,3),\:\: (2,3,4),\:\: (2,3,5),
\end{displaymath}
(i.e. corresponding to spherical triangular groups), then we remove the vertex.
Let $\hat{P}$ denote the differentiable orbifold obtained from $P$ with faces silvered, edge orders $n_e$, and with vertices removed as above.
We say that $\hat{P}$ has a \emph{Coxeter orbifold structure}.
For example, let $P$ be a convex hyperbolic polyhedron with dihedral angles
submultiples of $\pi$; we call $P$ a \emph{Coxeter polyhedron}.
Then $P$ will naturally have a Coxeter orbifold structure $\hat P$.

Now let $\mathfrak{D}(\hat{P})$ denote the deformation space of real projective structures on the Coxeter $3$-orbifold $\hat{P}$.
The work of Vinberg \cite{Vinberg1971} implies that each element of $\mathfrak{D}(\hat P)$ gives a convex projective structure (see Theorem 2 of \cite{Choi2006}).
That is, the image of the developing map of the orbifold universal cover
of $\hat P$
is projectively isomorphic to a convex domain in $\mathbb{R} P^3$
and the holonomy is a discrete faithful representation.
(For a precise definition of convexity, see \cite{Choi2006}.)

A point $p$ of $\mathfrak{D}(\hat{P})$ gives a fundamental polyhedron $P$ in $\R P^3$,
well defined up to projective automorphisms. We concentrate on the space
of $p \in \mathfrak{D}(\hat{P})$ giving
a {\em fixed} fundamental polyhedron $P$.
This space is called the \emph{restricted deformation space} of $\hat{P}$ and denoted by $\mathfrak{D}_P(\hat{P)}$. A point $t$ in $\mathfrak{D}_P(\hat{P)}$ is said to be \emph{hyperbolic} if it is given by a hyperbolic structure on $\hat{P}$.

\begin{definition}\label{def:deformable}
Let $P$ be a 3-dimensional hyperbolic Coxeter polyhedron, and let $\hat{P}$ denote its Coxeter orbifold structure.
Suppose that $t$ is the corresponding hyperbolic point of $\mathfrak{D}_P(\hat{P})$.
We call a neighbourhood of $t$ in $\mathfrak{D}_P(\hat{P})$ the \emph{local restricted deformation space} of $P$. We say that $\hat{P}$ is \emph{projectively deformable relative to the mirrors}, or
simply \emph{deforms rel mirrors},
if the dimension of its local restricted deformation space is positive.
Conversely, we say that $\hat{P}$ is \emph{projectively rigid relative to the mirrors}, or  \emph{rigid rel  mirrors}, if the dimension of
its local restricted deformation space is $0$.
\end{definition}

Choi  \cite{Choi2006} found a class of Coxeter 3-orbifolds whose restricted deformation spaces are understandable:
the {orderable} Coxeter orbifolds of normal type.
A Coxeter orbifold $\hat{P}$ is said to be  \emph{orderable} if the faces of $P$ can be ordered
so that each face contains
at most three edges that are edges of order 2 or edges in a face of higher index.
(See \S \ref{subsec:orderable} for the details, and for the definition of normal type.)

In this paper, we will study Coxeter orbifolds that are not orderable.
The following theorem describes the local restricted deformation space for a class of
Coxeter orbifolds arising from  {\em ideal} hyperbolic polyhedra, i.e. polyhedra with all vertices on the sphere at infinity.

\begin{theorem}\label{thm:main}
Let $P$ be an ideal $3$-dimensional hyperbolic polyhedron
whose dihedral angles are all equal to $\pi/3$, and suppose that $\hat{P}$ is given its Coxeter orbifold structure.
If $P$ is not a tetrahedron,
then a neighbourhood of the hyperbolic point in
$\mathfrak{D}_P(\hat{P)}$ is a smooth $6$-dimensional manifold.
\end{theorem}

If $P$ is a regular ideal tetrahedron then Theorem 3 of \cite{Choi2006}, obtained by J.R.~Kim in his master's thesis,
shows that $\mathfrak{D}_P(\hat{P)}$ is a $3$-dimensional cell.

The main ideas in the proof of Theorem \ref{thm:main} are as follows.
We first show that
$\mathfrak{D}_P(\hat{P)}$
is isomorphic to the solution set of  a system of polynomial equations following ideas of Vinberg \cite{Vinberg1971} and Choi \cite{Choi2006}. Since the faces of $P$ are fixed,
each projective reflection
in a face of the polyhedron is determined by a reflection vector $b_i$.
We then compute the Jacobian matrix of the equations for the $b_i$ at the hyperbolic point. This
reveals that the matrix has exactly the same rank as the Jacobian matrix of the equations for
the Lorentzian unit normals of a hyperbolic polyhedron with the given dihedral angles. By infinitesimal rigidity of the hyperbolic structure on $\hat P$, this matrix is of full rank and has kernel of dimension six; the result then follows from the  implicit function theorem.
 In fact, we can interpret the infinitesimal projective deformations as applying infinitesimal hyperbolic isometries to the reflection vectors.

The other two main results of this paper use various theoretical and computational methods
to determine the local restricted deformation spaces of Coxeter orbifolds arising from certain compact cubes and dodecahedra in hyperbolic $3$-space. These cubes and dodecahedra were chosen since they are workable using our methods, but not trivially.
The results are summarized in the following two theorems; the details are given in \S \ref{subsec:cubes} and \S \ref{subsec:dodec} below.

\begin{theorem}\label{thm:cube}
Consider the compact hyperbolic cubes such that each dihedral angle is $\pi/2$ or $\pi/3$.
Up to symmetries, there exist 34 cubes satisfying this condition.
For the corresponding hyperbolic Coxeter orbifolds,
10 are projectively deformable relative to the mirrors and the remaining 24 are projectively rigid relative to the mirrors.
\end{theorem}

\begin{theorem}\label{thm:dod}
Consider the compact hyperbolic dodecahedra such that each dihedral angle is $\pi/2$ or $\pi/3$, and
each face
has at most two dihedral angles equal to $\pi/2$.
Up to symmetries, there exist 13 dodecahedra satisfying these conditions.
For the corresponding hyperbolic Coxeter orbifolds,
only 1 is projectively deformable relative to the mirrors and the remaining 12 are projectively rigid relative to the mirrors.
\end{theorem}

If a face has more than two
edges of order two, then the corresponding reflection
is determined.
For dodecahedra, we assumed this condition fails
for every face and tabulated the results. Without this restriction,
the list of the possible dodecahedra would become very large and many of these would
prove to be projectively rigid relative to the mirrors by the linear test presented later in \S \ref{subsec:lin_test}.
It is future work to complete the task of fully classifying the Coxeter orbifold structures on dodecahedra
and cubes that are projectively deformable relative to the mirrors.

To obtain Theorems \ref{thm:cube} and \ref{thm:dod}, the polyhedra were first enumerated by
using a Matlab program to check the conditions of Andreev's theorem.
(See cu.m and do.m in \cite{Lee}.)
The remaining computations were done by Mathematica.
In the case of cubes, we used exact algebraic computations. However, in the case of dodecahedra numerical computations were used.
The detailed results of computations by Mathematica can be found at the web page \cite{Lee}.

The remainder of this paper is organized as follows.

Section \ref{prelim} reviews some well-known facts.
In \S \ref{subsec:orient_proj} we introduce {oriented projective structures} which
are in one-to-one correspondence with real projective structures.
In \S \ref{subsec:Vinberg} we describe Vinberg's results giving the general conditions satisfied by $n$-dimensional real projective reflection groups.
In \S \ref{subsec:andreev} we recall Andreev's theorem characterizing the 3-dimensional hyperbolic polyhedra
of finite volume with dihedral angles at most $ \pi/2$.
In \S \ref{subsec:orderable} we describe the results of Choi \cite{Choi2006} on the restricted deformation spaces of 3-dimensional Coxeter orbifolds that are orderable.

Section \ref{sec:restrict_def_space} identifies the restricted deformation space of real projective structures on a Coxeter orbifold $\hat P$ with the solution space of a system of polynomial equations.
In  \S  \ref{subsec:restict_rep_space} we introduce a space of restricted representations of the
orbifold fundamental group $\pi_1(\hat P)$.
In \S  \ref{subsec:proj_rep_space} we show that this representation space can be identified with the solution space of some polynomial equations as given by Vinberg.
In \S  \ref{subsec:def_vs_hol} we prove that the restricted deformation space
can identified with the set of solutions of Vinberg's equations, when the underlying
convex polyhedron $P$ has a discrete projective automorphism group.
In \S \ref{subsec:hyp_point} we recall the description of convex hyperbolic polyhedra by their Gram matrices, and use this to identify the
solutions to Vinberg's equations corresponding to a hyperbolic structure.

Section \ref{loc_defs} discusses general facts concerning a neighbourhood of a hyperbolic structure in the restricted deformation space of real projective structures on a Coxeter $3$-orbifold. The results from  \S \ref{subsec:def_vs_hol} show that this restricted deformation space is the solution space of a system of polynomial equations.
In  \S \ref{subsec:Z_tang_space} we study the Zariski tangent space of this solution space, and prove some general results on
local restricted deformation spaces.
In \S \ref{Z_tang_hyperbolic} we study the Zariski tangent space for the equations defining a hyperbolic structure.
In \S \ref{subsec:Thm1_proof} we compare the two Zariski tangent spaces at a hyperbolic point, and use
Garland-Raghunathan-Weil infinitesimal rigidity (\cite{GarlandRag}, \cite{Weil1962})
to prove Theorem \ref{thm:main}. (See Kapovich \cite{Kapovich1994} for a similar work in the conformally flat structures.)
In \S \ref{subsec:prisms}  we construct families of compact hyperbolic prism orbifolds, with number of faces arbitrarily large,
that are deformable relative to the mirrors but non-orderable.  In contrast, we note that orderable 3-dimensional compact hyperbolic
Coxeter polyhedra are always projectively rigid relative
to the mirrors when the number of faces is greater than $7$.

Section \ref{computations} is concerned with computing the dimension of local restricted
deformation spaces for hyperbolic Coxeter orbifolds corresponding to cubes and dodecahedra.
We carry out most of the computations using Mathematica. An outline of the computational algorithm is given in \S \ref{subsec:main_alg}.
In \S \ref{subsec:lin_test} we provide a simple test for the projective rigidity rel mirrors of 3-dimensional Coxeter polyhedra in real projective space.
In \S \ref{subsec:notation} we describe the notation used in figures and tables in this paper.
In \S \ref{subsec:cubes}-\ref{subsec:dodec} we give details of the methods used, and provide
detailed tables listing the dimensions of local restricted deformation spaces of cubes and dodecahedra.
The results show that computation of the Zariski tangent space is often sufficient; but  in other cases,
Gr\"obner bases are used to determine the structure of the local restricted deformation spaces.

This work was carried out while two of the authors were visiting the Department of Mathematics and Statistics of
the University of Melbourne; one of the authors also visited the Department of Mathematical and Computing Sciences in
Tokyo Institute of Technology. We thank both institutions very much for their kind hospitality towards the authors. We have benefited much from discussions with Yves Benoist, Bill Goldman, Misha Kapovich, who studied this kind of question before,
Sadayoshi Kojima, Daryl Cooper, and Steve Kerckhoff. We also very much appreciate the help and support
of many mathematicians and others whose names are not mentioned here.

\section{Preliminaries}\label{prelim}

This section gives the basic background material used in this article.

In \S  \ref{subsec:orient_proj} we give an alternative description of real projective structures that will be more convenient for us in this article, based on the projective sphere ${\mathbb S}^n$ and its group of projective  transformations $SL_\pm(n+1,\mathbb{R})$.
In \S \ref{subsec:Vinberg} we describe Vinberg's results giving the conditions
under which an $n$-dimensional Coxeter orbifold $\hat P$ admits a real projective structure.
This is equivalent to showing that the (orbifold) fundamental group $\pi_1(\hat{P})$ is isomorphic to a discrete subgroup of $SL_\pm(n+1,\mathbb{R})$.
We then concentrate on the case where $\hat P$ is a 3-dimensional
Coxeter orbifold. In \S \ref{subsec:andreev} we recall Andreev's theorem which explains
when $\hat P$ admits a finite volume hyperbolic structure.
Finally, in \S \ref{subsec:orderable}, Choi's results on the restricted deformation spaces of real projective structures on
\emph{orderable} Coxeter 3-orbifolds are described.

\subsection{Oriented real projective geometry}\label{subsec:orient_proj}
Instead of working in the $n$-dimensional real projective space
$\mathbb{R}P^n$,
it will be more convenient for us to work in the projective sphere $\mathbb{S}^n$, i.e. the set of
rays through the origin in $\R^{n+1}$.
As a $(G,X)$-structure, an {\em oriented projective structure} is a
$(SL_{\pm}(n+1,\mathbb{R}),\mathbb{S}^n)$-structure,
where
$$SL_{\pm}(n+1,\mathbb{R})=\{A \in GL(n+1,\R) : \det A = \pm 1\}$$
is the group of projective transformations of $\mathbb{S}^n$.
Recall that $\mathbb{S}^n$ double covers $\mathbb{R}P^n$ and
$SL_{\pm}(n+1,\mathbb{R})$ double covers $PGL(n+1,\mathbb{R})$.
A projective structure on an orbifold corresponds to a unique oriented projective structure and vice versa
(see \cite{Choi2006} and \cite[p. 143]{Thurston_book}).
From now on, by a real projective structure, we always mean an oriented projective structure.

There is a one-to-one correspondence between the set of vector subspaces of $\R^{n+1}$
and the set of great spheres in $\mathbb{S}^n$.
In particular, a 1-dimensional subspace corresponds to a pair of antipodal points
and an $n$-dimensional subspace gives a great $(n-1)$-sphere in $\mathbb{S}^n$.
Further,  a component of the complement of a great $(n-1)$-sphere (i.e. an open hemisphere)
can be identified with an affine $n$-space. We call this an {\em affine patch}.

In this paper, we define a {\em convex polytope} $P$ in $\mathbb{S}^n$ to be a
precompact convex polytope in an affine patch of $\mathbb{S}^n$. The image of
such a polytope under the double covering is called a convex polytope in $\mathbb{R}P^n$.
We define $k(P)$ as the dimension of the subgroup of $SL_{\pm}(n+1,\mathbb{R})$
preserving $P$. This is the same as the dimension of the group of projective automorphisms of
the image of $P$ under the double-covering map.

Hyperbolic geometry arises naturally as a sub-geometry of oriented projective geometry.
Let $\langle \cdot \, , \cdot \rangle$ denote the Lorentzian inner product on $\R^{n+1}$ defined by
 \begin{displaymath}
 \langle x, y\rangle = -x_1 y_1 + x_2 y_2 +  \ldots + x_{n+1} y_{n+1},
 \end{displaymath}
 and let $B' \subset \mathbb{S}^n$ be the open $n$-ball consisting of
 rays through the origin in the cone $\{ x \in \R^{n+1} : \langle x, x \rangle < 0, ~~x_1 > 0 \}$.
  Then we can regard hyperbolic space $\mathbb{H}^n$ as the open ball $B'$, and
 the group of hyperbolic isometries $\text{Isom}(\H^n)$ is the subgroup $O_0(1,n)$
 of $SL_{\pm}(n+1,\R)$ preserving $B'$.

Radial projection maps
$B'$ diffeomorphically to an open $n$-ball $B$
in the affine hyperplane $x_1=1$, and $\text{Isom}(\H^n)$ corresponds to the
closed subgroup $PO(1, n) \subset PGL(n+1, \mathbb{R})$ of projective automorphisms of $B$.
This gives the Klein model for hyperbolic geometry.

Alternatively, hyperbolic space $\mathbb{H}^n$ can be embedded in
$\mathbb{R}^{n+1}$ as the upper sheet of a hyperboloid
$$
\langle x, x \rangle=-1, ~~x_1>0,
$$
and $\text{Isom}(\H^n)$ is the subgroup $O_0(1,n)$ of $SL_{\pm}(n+1,\R)$ preserving $\H^n$.

\subsection{Vinberg's results}\label{subsec:Vinberg}

This subsection gives a summary of results from Vinberg's article \cite{Vinberg1971}.
An alternative treatment is given in Benoist's notes \cite{BenoistNote}.
Vinberg gave the general conditions under which a
Coxeter orbifold admits a real projective structure, and a criterion to decide whether it is a hyperbolic structure or not.

Let $V$ be the $(n+1)$-dimensional real vector space $\mathbb{R}^{n+1}$.
A (projective) \emph{reflection} $R$ is an element of order $2$ of
$SL_{\pm}(n+1,\mathbb{R})$
which is the identity on a hyperplane $U$. All reflections are of the form
\begin{displaymath}
R=Id-\alpha \otimes b
\end{displaymath}
for some linear functional $\alpha \in V^*$ and a vector $b \in V$ with $\alpha(b)=2$.
Here, the kernel of $\alpha$ is the subspace $U$ of fixed points of $R$ and $b$ is
the \emph{reflection vector}, i.e. an eigenvector corresponding to the eigenvalue $-1$.
A \emph{rotation} is an element of $SL_{\pm}(n+1,\mathbb{R})$ which is the identity on
a subspace of codimension 2 and is conjugate to a matrix
\begin{displaymath}
 \left( \begin{array}{cc}
\cos\theta & -\sin\theta \\
\sin\theta & \cos\theta \end{array}\right)
\end{displaymath}
in a suitable supplementary basis. The real number $\theta \in [0,\pi]$ is the \emph{angle} of the rotation.

We consider $\mathbb{S}^n$ as the set of rays in  $\mathbb{R}^{n+1} $ from the origin.
Let $P$ be an $n$-dimensional convex polytope in $\mathbb{S}^n$ and for each (codimension one) face $F_i$ of $P$, take
a linear functional $\alpha_i$ for $F_i$ and choose a projective reflection $R_i=Id-\alpha_i \otimes b_i$ with $\alpha_i(b_i)=2$ which fixes $F_i$.
By making a suitable choice of signs, we will assume that $P$ is defined by the inequalities
\begin{equation}
\alpha_i \leq 0 \hspace{10mm} i=1, \ldots, f,
\end{equation}
where $f$ is the number of faces of $P$.
The group $\Gamma \subset SL_{\pm}(n+1,\mathbb{R})$ generated by all these reflections $R_i$ is called a \emph{linear Coxeter group} if
$$ \gamma P^\circ \cap P^\circ = \varnothing \:\;\:\text{for every}\: \gamma \in \Gamma \backslash \{1\},$$
where $P^\circ$ is the interior of $P$.
The $f \times f$ matrix $A=(a_{ij})$, $a_{ij}=\alpha_i(b_j)$, is called the \emph{Cartan matrix} of $\Gamma$. Vinberg proved that the following conditions are necessary and sufficient for $\Gamma$ to be a linear Coxeter group:
\begin{enumerate}
\item[(C1)] $a_{ij} \leq 0$ for $i \neq j$, and $a_{ij}=0 \Leftrightarrow a_{ji}=0$.
\item[(C2)] $a_{ii}=2$; and for $i\ne j$, $a_{ij}a_{ji}\geq 4$ or $a_{ij}a_{ji}=4\cos^2\frac{\pi}{n_{ij}}$, $n_{ij}$ an integer.
\end{enumerate}
In fact, if $a_{ij}a_{ji}=4\cos^2\frac{\pi}{n_{ij}}$ then the product $R_i R_j$ is a rotation of angle $2\pi/n_{ij}$ and the group generated by two reflections $R_i$ and $R_j$ is the dihedral group $D_{n_{ij}}$. Note that (C1) and (C2) imply that $a_{ij}=a_{ji}=0$ if $n_{ij}=2$; however $a_{ij}\ne a_{ji}$ in general when $n_{ij}>2$.

For each reflection $R_i$,  $\alpha_i$ and $b_i$ are defined up to transformations
\begin{displaymath}
\alpha_i \mapsto d_i\alpha_i \;\; \text{and} \;\; b_i \mapsto d_i^{-1}b_i  \text{ with } d_i  > 0.
\end{displaymath}
Hence the Cartan matrix of $\Gamma$ is defined up to conjugation by a diagonal matrix
with positive diagonal entries.

For any $x \in P$, let $\Gamma_x$ denote the subgroup of $\Gamma$ generated by reflections in those faces of $P$ which contain $x$. Define $P^f=\{x \in P \,|\, \Gamma_x \text{ is finite} \}$. Then the following statements are true:
\begin{enumerate}
\item $C=\cup_{\gamma \in \Gamma} \gamma P$ is convex.
\item $\Gamma$ is a discrete subgroup of $SL_{\pm}(n+1,\mathbb{R})$ preserving $C^\circ$.
\item $C^\circ \cap P = P^f$, and is homeomorphic to $C^\circ/ \Gamma$.
\end{enumerate}

 Thus $C^\circ$ gives a convex open subset of the projective sphere $\mathbb{S}^n$, and
 $C^\circ/ \Gamma$ determines a convex real projective structure on the Coxeter orbifold
 $\hat{P}$ associated with $P$.

To state the next theorem, we introduce the following notation and definitions:
if $X=(X_1,\ldots, X_{n+1}) \in V$,
we write $X>0$ if $X_i>0$ for every $i$, and $X \geq 0$ if $X_i \geq 0$ for every $i$.
A matrix $A$ is of \emph{negative type} if there exists $X >0$ such that $AX<0$, and if $X \geq 0$ and $AX \geq 0$ imply $X=0$.
A matrix $A$ is \emph{indecomposable} if it cannot be represented as a direct sum of two matrices.
Two matrices $A$ and $B$ are said to be \emph{equivalent} if $A=DBD^{-1}$ for a diagonal matrix $D$ having positive entries.
A linear Coxeter group $\Gamma$ is called a \emph{hyperbolic} Coxeter group
if $\Gamma$ is derived from a discrete group generated by reflections in $\mathbb{H}^n$, and no proper plane of $\mathbb{H}^n$ or any point at infinity is $\Gamma$-invariant.

\begin{theorem} \label{thm:vinberg} {\rm (Vinberg \cite{Vinberg1971})}
A linear Coxeter group $\Gamma$ is hyperbolic if and only if
the Cartan matrix $A$ of $\Gamma$ is indecomposable, of negative type, and equivalent to a symmetric matrix of signature $(1,n)$.
\end{theorem}

\subsection{Andreev's theorem}\label{subsec:andreev}
The 3-dimensional Coxeter orbifolds which admit a finite volume hyperbolic structure
have been classified by Andreev \cite{Andreev1970,Andreev19702}.

Let $X$ be an $3$-dimensional space of constant curvature, with  group of isometries denoted Isom($X$). A convex polyhedron $P$ in $X$ is called a \emph{Coxeter polyhedron} if all the dihedral angles of $P$ are submultiples of $\pi$. Let $P$ be a Coxeter polyhedron, and $\Gamma$ be the group generated by reflections in its faces. Then $\Gamma$ is a discrete subgroup of $\text{Isom}(X)$, and $P$ is its fundamental polyhedron. Conversely, every discrete subgroup of $\text{Isom}(X)$ generated by reflections can be obtained in this manner.

A nice property of a Coxeter polyhedron is that its dihedral angles are \emph{non-obtuse}, i.e. the dihedral angles do not exceed $\pi/2$.
In 1970, E.M. Andreev \cite{Andreev1970} gave a full description of $3$-dimensional compact hyperbolic polyhedra with non-obtuse dihedral angles.

Let $C$ be an {abstract} $3$-dimensional polyhedron and $C^*$ be its dual.
A simple closed curve $\gamma$ is called a \emph{$k$-circuit} if it consists of $k$ edges of $C^*$.
A circuit $\gamma$ is \emph{prismatic} if all of the endpoints of the edges of $C$ which $\gamma$ meets are different.

Suppose that $C$ is \emph{not a tetrahedron} and non-obtuse angles $\theta_{ij} \in (0,\pi/2]$ are given corresponding to each edge $F_{ij}=F_i \cap F_j$ of $C$, where $F_i$ are the faces of $C$. Then the following conditions $(A1)$--$(A4)$ are necessary and sufficient for the existence of a compact $3$-dimensional hyperbolic polyhedron $P$ which realizes $C$ with dihedral angle $\theta_{ij}$ at each edge $F_{ij}$.
\begin{enumerate}
\item[(A1)] If $F_{ijk}=F_i \cap F_j \cap F_k$ is a vertex of $C$ then
    $$ \theta_{ij} + \theta_{jk} + \theta_{ki} > \pi.$$
\item[(A2)] If $F_i$, $F_j$, $F_k$ form a prismatic $3$-circuit, then
    $$ \theta_{ij} + \theta_{jk} + \theta_{ki} < \pi.$$
\item[(A3)] If $F_i$, $F_j$, $F_k$, $F_l$ form a prismatic $4$-circuit, then
    $$ \theta_{ij} + \theta_{jk} + \theta_{kl} + \theta_{li} < 2\pi.$$
\item[(A4)] If $C$ is a triangular prism with
triangular faces $F_1$ and $F_2$, then
        $$ \theta_{13} + \theta_{14} + \theta_{15} + \theta_{23} + \theta_{24} + \theta_{25} < 3\pi. $$
\end{enumerate}
Furthermore, this polyhedron is unique up to hyperbolic isometries.

\medskip
Andreev \cite{Andreev19702} also showed that the following conditions $(\widetilde{A1})$--$(\widetilde{A6})$ are necessary and sufficient for the existence of a $3$-dimensional hyperbolic polyhedron $P$ of finite volume which realizes $C$ with dihedral angle $\theta_{ij} \in (0,\pi/2]$ at each edge $F_{ij}$.
\begin{enumerate}
\item[$(\widetilde{A1})$] If $F_{ijk}=F_i \cap F_j \cap F_k$ is a vertex of $C$ then
    $$ \theta_{ij} + \theta_{jk} + \theta_{ki} \geq \pi.$$
\item[$(\widetilde{A2})$] If $F_{ijkl}=F_i \cap F_j \cap F_k \cap F_l$ is a vertex of $C$ then
    $$ \theta_{ij}+\theta_{jk}+\theta_{kl}+\theta_{li}=2\pi.$$
\item[$(\widetilde{A3})$] is the same as (A2).
\item[$(\widetilde{A4})$] is the same as (A3).
\item[$(\widetilde{A5})$] is same as (A4).
\item[$(\widetilde{A6})$] If $F_i$, $F_j$, $F_k$ are faces with $F_i$ and $F_j$ adjacent, $F_j$ and $F_k$ adjacent, and $F_i$ and $F_k$ are not adjacent
but meet in a vertex not in $F_j$, then
    $$ \theta_{ij}+\theta_{jk} < \pi.$$
\end{enumerate}
Again, the hyperbolic polyhedron is unique up to hyperbolic isometries.

\medskip
Note that if the vertices of $C$ are all trivalent then conditions $(\widetilde{A2})$ and $(\widetilde{A6})$ are not needed.

\subsection{Orderability results}\label{subsec:orderable}

This subsection describes the main theorem of \cite{Choi2006}. As we mentioned in the introduction, if a Coxeter orbifold $\hat{P}$ satisfies the condition of orderability, then we understand the restricted deformation space of real projective structures on $\hat{P}$.

Let $P$ be a fixed $3$-dimensional convex polyhedron in $\mathbb{S}^3$
with given edge orders, and let $\hat{P}$ be the corresponding Coxeter orbifold. Denote the numbers of vertices, edges and faces of $P$ by $v$, $e$, $f$ respectively. Let $e_2$ be the number of edges of order 2 in $\hat P$. Let $k(P)$ be the dimension of the group of projective automorphisms of $P$. Then $k(P)= 3$ if $P$ is tetrahedron,
$k(P) =1$ if $P$ is the cone over a polygon other than a triangle, and $k(P)=0$ otherwise. (See Lemma \ref{lem:kp} below.)

The orbifold $\hat{P}$ is called a \emph{normal-type} Coxeter orbifold if it is
{\bf not} one of the following types:
\begin{itemize}
\item  a  \emph{cone-type Coxeter orbifold}, whose underlying polyhedron is topologically a cone from a face  $F$ to a vertex, and all edges of $F$ have edge order $2$,
\item a \emph{product-type Coxeter orbifold}, whose underlying polyhedron is topologically a polygon times an interval, and all edges of the top and bottom faces have edge order $2$,
\item a Coxeter orbifold with finite fundamental group.
\end{itemize}

Recall that a Coxeter orbifold $\hat{P}$ is said to be \emph{orderable} if the faces of $P$ can be ordered
so that each face contains
at most three edges that are edges of order 2 or edges in a face of higher index.
Then we have:

\begin{theorem}\label{thm:choi} {\rm (\cite{Choi2006}) }
Let $P$ be a 3-dimensional convex polyhedron in $\mathbb{S}^3$ and let $\hat{P}$ be given a normal-type Coxeter orbifold structure. Suppose that $\hat{P}$ is orderable. Then the restricted deformation space of projective structures on $\hat{P}$ is
a smooth manifold of dimension $3f-e-e_2-k(P)$ if it is not empty.
\end{theorem}

Examples of orderable orbifolds are obtained if $P$ is any convex polyhedron with all faces triangular.
An antiprism (i.e. drum-shaped convex polyhedron
with $n$-gons on the top and bottom joined up by a band of $2n$-triangles)
with arbitrary orders given to the edges is orderable,
since we can order the top and the bottom faces to have the highest two indices.
By Andreev's theorem, an antiprism with all angles $\pi/2$ admits a complete hyperbolic Coxeter orbifold structure (see also Thurston \cite{Thurston}).
A triangular prism carries compact hyperbolic Coxeter orbifold structures
and these are all orderable.

However the cube and dodecahedron do not carry an orderable Coxeter orbifold structure, since a lowest index face in an orderable orbifold must be triangular.

\section{The restricted deformation space of real projective structures}\label{sec:restrict_def_space}

In this section, the restricted deformation space of real projective structures on an $n$-dimensional Coxeter orbifold $\hat P$ is discussed, and identified with a space of representations.

 In \S  \ref{subsec:restict_rep_space} we define a suitable space of restricted representations from
 $\pi_1(\hat P)$ into $SL_\pm(n+1,\R)$.
  In \S  \ref{subsec:proj_rep_space} we show that this restricted representation space can be identified with the solution space of a system of polynomial equations given by Vinberg (Proposition \ref{prop:eqnsol}).
In \S  \ref{subsec:def_vs_hol} we prove that the restricted deformation space is homeomorphic to the space of restricted representations and to the set of solutions of Vinberg's equations (Theorems \ref{isotopy_equiv} and   \ref{thm:fundamental}), when the underlying convex polyhedron $P$ has a discrete projective automorphism group.
In \S  \ref{subsec:hyp_point} we look at the equations satisfied by the Lorentzian unit normals to a hyperbolic polyhedron, and show that a hyperbolic structure on $\hat P$ corresponds to the single point
in the solution space of \S  \ref{subsec:proj_rep_space}.

\subsection{The restricted representation space}
\label{subsec:restict_rep_space}

Let $P$ be a fixed $n$-dimensional convex polytope contained in $\mathbb{S}^n$,
and $\hat{P}$ an associated Coxeter orbifold.
We now
identify the deformation space $\mathfrak{D}(\hat{P})$ of real projective structures on $\hat{P}$ with the deformation space of $(G,X)$-structures on $\hat{P}$, where
$G =SL_{\pm}(n+1,\mathbb{R})$ and $X =\mathbb{S}^n$
is the projective sphere in $V=\R^{n+1}$.

Sending a development pair $(D,h)$ to its holonomy representation $h$ induces a local homeomorphism
$$
hol\: : \: \tilde{\mathfrak{D}}(\hat{P}) \rightarrow
\text{Hom}(\pi_1(\hat{P}),G),
$$
where  $\tilde{\mathfrak{D}}(\hat{P})$ denotes the space of isotopy-equivalence classes of development pairs.
(See Theorem 1 of \cite{Choi2004} and Proposition 3 of \cite{Choi2006}, where only a sketch proof is given.
One of the authors is writing a more complete account in a generalized setting \cite{Choi2010}.)

Recall that $\pi_1(\hat P)$ is a Coxeter group with standard generator $r_i$
corresponding to the $i$th face of $P$. To study the restricted deformation space
we consider the subset
$$\tilde{\mathfrak{D}}_{P}(\hat{P}) \subset \tilde{\mathfrak{D}}(\hat{P})$$
giving projective structures {\em with fundamental polyhedron $P$}. More precisely, let $H_i$ denote
the hyperplane in $V$ containing the $i$th face of $P$.
Then $\tilde{\mathfrak{D}}_{P}(\hat{P})$ consists of the isotopy-equivalence classes
$[(D,h)]$ of developing pairs $(D,h)$ such that each $h(r_i)$ is a reflection with fixed point set
$\text{Fix}(h(r_i))= H_i$.

\begin{lemma}
For every $[(D,h)] \in \tilde{\mathfrak{D}}_{P}(\hat{P})$
the holonomy $h$ lies in the subset
$$\text{Hom}_P(\pi_1(\hat{P}),G) \subset
\text{Hom}(\pi_1(\hat{P}),G)$$
consisting of representations such that each $h(r_i)$ is a projective reflection fixing $H_i$, and $h(r_i r_j)$ is a rotation by $2\pi /n_{ij}$ whenever $F_i \cap F_j$ is a codimension $2$ face of $P$ of order $n_{ij}$.
\end{lemma}

\begin{proof} The definition of the orbifold structure on $\hat P$ shows that
the local action of $r_i$ and $r_j$ on the universal cover of $\hat P$ is given by a standard dihedral group of order $2 n_{ij}$, generated by involutions fixing two hypersurfaces meeting transversally at an angle $\pi/n_{ij}$. Given a real projective structure on $\hat P$, this action is transferred by the developing map into
$\mathbb{S}^n$. Hence each $h(r_i)$ is a projective reflection, and $h(r_i)h(r_j)$ is conjugate to a rotation by $\pi/n_{ij}$.
\end{proof}

We call $\text{Hom}_P(\pi_1(\hat{P}),G)$ the space of \emph{restricted representations} from $\pi_1(\hat P)$ to $G$. Lemma 1 shows that $hol$ restricts to a map
$$hol_P : \tilde{\mathfrak{D}}_{P}(\hat{P}) \to \text{Hom}_P(\pi_1(\hat{P}),G). $$
In Theorem \ref{isotopy_equiv}, we will show that this is a homeomorphism.

\subsection{Restricted representations and Vinberg's equations}\label{subsec:proj_rep_space}

We now give a very explicit description of the restricted representation space
 $\text{Hom}_P(\pi_1(\hat{P}),G)$ for the fundamental group of a
 Coxeter orbifold $\hat P$.
Let $V=\R^{n+1}$
and let $P$ be a \emph{fixed} convex polytope in $\mathbb{S}^n$.

Assume that $P$ is given by a system of linear inequalities
\begin{displaymath}
 \alpha_i \leq 0, \hspace{3mm} i=1,\ldots, f,
\end{displaymath}
where $\alpha_i \in V^*$ and $f$ is the number of codimension one faces of $P$.
Suppose $b_i \in V$
for $1 \leq i \leq f$ are reflection vectors with $\alpha_i(b_i)=2$.
Let $R_i$ be the reflections defined by $R_i=Id-\alpha_i \otimes b_i$ for $i=1,...,f$,
and let $\Gamma \subset SL_{\pm}(n+1,\mathbb{R})$ be the group generated by the $R_i$. Then the matrix $A = (a_{ij}) = (\alpha_i(b_j))$ is the $f \times f$ Cartan matrix of $\Gamma$.

We  now fix  orders $n_{ij}$ for the codimension 2 faces of $P$ consider the \emph{restricted} deformation space of the corresponding Coxeter orbifold $\hat P$.
Then the $\alpha_i$'s will be \emph{fixed}, and $b_i$'s are \emph{variables,}
so Vinberg's result leads us to solve the following system of polynomial equations:
\begin{itemize}
\item For each $i=1,\ldots,f$,
\begin{equation}\label{eqn:normalize}
a_{ii} = \alpha_i(b_i) = 2,
\end{equation}
\item If $F_i$ and $F_j$ are adjacent in $P$ and $n_{ij}>2$,
\begin{equation}\label{eqn:Vinberg1}
a_{ij}a_{ji}= \alpha_i(b_j)\alpha_j(b_i)=4\cos^2(\pi/n_{ij}),
\end{equation}
\item If $F_i$ and $F_j$ are adjacent in $P$ and $n_{ij}=2$,
\begin{equation}\label{eqn:Vinberg2}
a_{ij}=\alpha_i(b_j)= 0 \;\; \text{and} \;\;  a_{ji}=\alpha_j(b_i)=0.
\end{equation}
\end{itemize}
(Note the difference between the cases $n_{ij}=2$ and $n_{ij}>2$.)

We call these polynomial equations (\ref{eqn:normalize})--(\ref{eqn:Vinberg2}) the \emph{Vinberg equations}.
Let $N$ be the number of Vinberg equations and let
$\Phi_{\hat P}\,:\, V^f=\mathbb{R}^{(n+1)f} \rightarrow \mathbb{R}^N$ be the map given by
\begin{displaymath}
(b_1, \ldots, b_f) \mapsto (\Phi_1, \ldots, \Phi_N),
\end{displaymath}
where $\{\Phi_k\}_{k=1}^{N}$ is the set of polynomials $a_{ii}-2$, $a_{ij}a_{ji}-4\cos^2(\pi/n_{ij})$,
or $a_{ij}$, $a_{ji}$ as in the above equations (\ref{eqn:normalize})--(\ref{eqn:Vinberg2}).
Note that $N = f + e + e_2$, where
$e$ is the number of codimension 2 faces of $P$, and $e_2$ is the number of codimension 2 faces of order 2.

\begin{proposition}\label{prop:eqnsol}
Suppose that the linear functionals $\alpha_i$ defining  the faces of $P$ are fixed.
Then there is a homeomorphism
$$
\mathcal{H}: \Phi^{-1}_{\hat P}(0)  \rightarrow  \text{Hom}_P(\pi_1(\hat{P}),SL_{\pm}(n+1,\mathbb{R}))
$$
where $\mathcal{H}$ sends
$(b_1,\ldots, b_f)$ to the homomorphism $h$ with $h(r_i)= Id - \alpha_i \otimes b_i$.
The map $\mathcal{H}$ is a polynomial map with a rational inverse $\mathcal{R}$.
\end{proposition}
\begin{proof}
Solving the Vinberg equations (\ref{eqn:normalize})--(\ref{eqn:Vinberg2}) is equivalent to finding
reflections $R_i$, $i=1,...,f$, corresponding to the faces of $P$,  such that
 $R_i R_j$ is conjugate to a rotation by $2\pi/n_{ij}$  whenever $F_i $ and $F_j$
 meet along a codimension 2 face. This follows, for example, from Lemma 1.2 of  \cite{BenoistNote}.

  Conversely, given a reflection matrix $R_i = Id - \alpha_i\otimes b_i$, the reflection
  vector $b_i$ is uniquely determined since $\alpha_i$ is fixed. In fact, $b_i$ is the unique
  eigenvector of $R_i$ with eigenvalue $-1$ satisfying the normalization condition
  $\alpha_i(b_i)=2$. It follows easily that the inverse map $\mathcal{R}$ taking each reflection matrix $h(r_i) = R_i$ to its reflection vector $b_i$ is a rational map.
  \end{proof}

From now on, the space of representations $\text{Hom}_P(\pi_1(\hat{P}),SL_{\pm}(n+1,\mathbb{R}))$ will be identified with $\Phi_{\hat P}^{-1}(0)$.

\subsection{The restricted deformation space}\label{subsec:def_vs_hol}

Let $P$ be a convex polytope in $\mathbb{S}^n$, and $\hat P$ an associated Coxeter orbifold.
In this section,  we will show that
that the restricted space of isotopy classes of real projective structures $\tilde{\mathfrak{D}}_{P}(\hat{P})$
can identified with the restricted representation space, and with the
set of solutions to Vinberg's equations (\ref{eqn:normalize})--(\ref{eqn:Vinberg2}).
In the generic case where the group of projective automorphisms of $P$ is discrete, these spaces
are also homeomorphic to the restricted deformation space ${\mathfrak{D}}_{P}(\hat{P})$ of
real projective structures on $\hat P$.

\begin{theorem}\label{isotopy_equiv}
The maps
 $$\tilde{\mathfrak{D}}_{P}(\hat{P}) \xrightarrow{hol_P} \text{Hom}_P(\pi_1(\hat{P}),G)
\xrightarrow{~\mathcal{R}~} \Phi_{\hat P}^{-1}(0)$$
 are homeomorphisms.
\end{theorem}

\begin{proof}
Given
a set of reflection vectors
$(b_1, \ldots, b_f) \in \Phi_{\hat P}^{-1}(0)$, the work of
Vinberg  (see \cite{Vinberg1971} or Theorem 1.5 of \cite{BenoistNote}) shows that
\begin{enumerate}
\item[(i)] the corresponding reflections
$R_i = Id -  \alpha_i\otimes b_i$
generate a discrete group $\Gamma$ isomorphic to $\pi_1(\hat P)$,
\item[(ii)] the images $\gamma P$ for $\gamma\in \Gamma$ tile an convex open subset $\Omega \subset \mathbb{S}^n$, and
\item[(iii)] the quotient orbifold $\Omega/\Gamma$ is isomorphic to $\hat P$.
\end{enumerate}
Thus we obtain a convex real projective structure on $\hat P$ and isotopy class of
development pair $[(D,h)] \in \tilde{\mathfrak{D}}_{P}(\hat{P})$
that maps to $(b_1, \ldots, b_f)$.
This gives continuous inverses to the maps in the theorem.
\end{proof}

Next we study the restricted deformation space $\mathfrak{D}_P(\hat{P)}$. Let
$$G_P = \{ g \in G : g(H_i)=H_i \text{ for all }i \}$$
be the subgroup of $G$
that preserves $P$ and each of its faces (and hence preserves each of its vertices).
Note that $\dim G_P = k(P)$, where $k(P)$ denotes the
dimension of the group of projective automorphisms of $P$  as in \cite{Choi2006}.

\begin{proposition}\label{restrict_deform}
The group  $G_P$ acts on $\tilde{\mathfrak{D}}_{P}(\hat{P})$,
and the quotient space
$\tilde{\mathfrak{D}}_{P}(\hat{P})/G_P$
 is homeomorphic to ${\mathfrak{D}}_{P}(\hat{P})$.
 \end{proposition}

 \begin{proof} We write  $\tilde{\mathfrak{D}} = \tilde{\mathfrak{D}}(\hat{P})$,
 $\tilde {\mathfrak{D}}_{P} = \tilde{\mathfrak{D}}_{P}(\hat{P})$, and
 $ {\mathfrak{D}}_{P} ={\mathfrak{D}}_{P}(\hat{P})$. Now $G$ acts on $\tilde{\mathfrak{D}}$
 by equation (\ref{PGLaction}), and we let
$\pi:   \tilde{\mathfrak{D}} \to  \tilde{\mathfrak{D}}/G = {\mathfrak{D}}$
be the natural quotient map.  Then
${\mathfrak{D}}_{P} =\pi(\tilde{\mathfrak{D}}_{P} ) = \hat{\mathfrak{D}}_{P}/G$,
 where $  \hat{\mathfrak{D}}_{P} =  \pi^{-1}({\mathfrak{D}}_{P})  = G \cdot \tilde {\mathfrak{D}}_{P} \subset
 \tilde {\mathfrak{D}}$.

Let $[(D,h)] \in \tilde {\mathfrak{D}}_{P}$
 and $g \in G$.  Then $g \cdot [(D,h)] = [(g \circ D, g \circ h \circ g^{-1})]$ and
 $\text{Fix}(g \circ h(r_i) \circ g^{-1}) = g \cdot \text{Fix}(h(r_i)) = g (H_i)$. Thus
 \begin{eqnarray} \label{GPaction}
g \cdot [(D,h)] \in \tilde{\mathfrak{D}}_{P} \text{ if and only if } g \in G_P .
\end{eqnarray}
In particular, it follows that  $G_P$ acts on $\tilde{\mathfrak{D}}_{P}$,
and we let $\pi_P:   \tilde{\mathfrak{D}}_P \to  \tilde{\mathfrak{D}}_P/G_P$
 be the natural quotient map.

Now the composition $\tilde{\mathfrak{D}}_P \subset \hat{\mathfrak{D}}_{P} \to \hat{\mathfrak{D}}_{P}/G$
is continuous and constant on $G_P$ orbits, so there is an induced
continuous map $\tilde{\mathfrak{D}}_P/G_P \to \hat{\mathfrak{D}}_{P}/G = \mathfrak{D}_{P}$,
taking $G_P\cdot y$ to $G \cdot y$ for $y \in \tilde{\mathfrak{D}}_P$. This is a bijection
by observation (\ref{GPaction}).  To finish the proof we show that the inverse is continuous.

First,
define $\phi : G \times \tilde{\mathfrak{D}}_{P} \to \hat{\mathfrak{D}}_{P}$ by $\phi(g,z)=g \cdot z$.
This is a continuous, surjective, open map, hence an identification map.
Let $p: G \times \tilde{\mathfrak{D}}_{P} \to \tilde{\mathfrak{D}}_{P}$ be the projection onto the second factor.
Now if $y,y' \in \phi^{-1}(x)$ then (\ref{GPaction}) shows that $p(y')=g \cdot p(y)$ for some $g \in G_P$, so $\pi_P \circ p(y)=
\pi_P \circ p(y')$. Hence $\pi_P \circ p \circ \phi^{-1}$ is well-defined and gives a continuous function
$ \hat{\mathfrak{D}}_{P} \to \tilde{\mathfrak{D}}_{P}/G_P$.
This is constant on $G$-orbits, so induces a continuous map
$\hat{\mathfrak{D}}_{P}/G \to \tilde{\mathfrak{D}}_{P}/G_P$. This is the desired inverse.
(More explicitly, the inverse
is given by $G\cdot x \mapsto (G\cdot x) \cap \tilde{\mathfrak{D}}_{P}$ for $x \in  \hat{\mathfrak{D}}_{P}$.)
  \end{proof}

Since the homeomorphism $hol_P$ in Theorem \ref{isotopy_equiv} is
equivariant with respect to the action of $G_P$, we also obtain the following.

\begin{corollary} The map $hol_P$ induces a homeomorphism
$${\mathfrak{D}}_{P}(\hat{P}) \cong \tilde{\mathfrak{D}}_{P}/G_P \to
 \text{Hom}_P(\pi_1(\hat{P}),G)/G_P ,$$
 where $G_P$ acts on $\text{Hom}_P(\pi_1(\hat{P}),G)$ by conjugation.
 \end{corollary}

In the remainder of this paper we concentrate on the generic case, where $k(P)=0$.
Then we have

\begin{theorem}\label{thm:fundamental}
If $k(P)=0$ then $G_P $ is a trivial group.  Hence we have homeomorphisms
$${\mathfrak{D}}_{P}(\hat{P})
\cong  \text{Hom}_P(\pi_1(\hat{P}),G) \cong \Phi_{\hat P}^{-1}(0).$$
\end{theorem}

\begin{proof}
If $k(P)=0$ then $G_P$ is a discrete group. Since $G_P$ acts trivially on each
face and vertex of $P$, it follows that the group is a trivial group.
The rest follows from Theorem \ref{isotopy_equiv}
and Proposition \ref{restrict_deform}.
\end{proof}

The following result shows that $k(P)=0$ holds for most convex polyhedra $P \subset \mathbb{S}^3$.

\begin{proposition}\label{lem:kp}
Let $P$ be a convex $3$-dimensional polyhedron in
$\mathbb{S}^3$.
Then
\begin{itemize}
\item $k(P) = 3$ if $P$ is a tetrahedron,
\item $k(P)=1$ if $P$ is a convex cone over a polygon other than a triangle, and
\item $k(P)=0$ otherwise.
\end{itemize}
Thus, $k(P)$ only depends on the
combinatorial type of $P$, and not on the geometric shape of $P$.
\end{proposition}

\begin{proof} It suffices to consider a projective
automorphism that fixes every vertex of $P$. If there is a face containing four or more
vertices, then every point of the face is fixed.  Such an automorphism is represented
by a diagonal matrix in $PGL(4,\R)$ with diagonal entries $1,1,1,\lambda$. Hence
$k(P)= 1$ if $P$ is a cone, and $k(P) =0$ otherwise.
If every face of $P$ is a triangle, then a similar argument shows that $k(P) = 3$ if $P$ is a tetrahedron, and $k(P) = 0$ otherwise.
The last sentence follows from these observations.

\end{proof}

\subsection{The hyperbolic point}\label{subsec:hyp_point}

Let $V$ be an $(n+1)$-dimensional vector space over $\mathbb{R}$ with coordinates $x_1, \ldots, x_{n+1}$, and let $P$ be a Coxeter polytope in Klein's model of $n$-dimensional hyperbolic space $\mathbb{H}^n$
with faces $F_i$ for $i=1, \ldots, f$.
Let $\nu_i \in V$ denote the outward unit normal to $F_i$
with respect to the Lorentzian inner product on $V$, defined by
 \begin{displaymath}
 \langle x, y\rangle = -x_1 y_1 + x_2 y_2 +  \ldots + x_{n+1} y_{n+1}.
 \end{displaymath}
 Then $P$ is defined by the system of linear inequalities
 \begin{displaymath}
 \langle \nu_i, x \rangle \leq 0, \hspace{3mm} i=1,\ldots, f, \;\;\text{and}\;\; x_1=1.
\end{displaymath}

Now the problem of constructing a hyperbolic polyhedron $P$ with prescribed dihedral angles $\pi/n_{ij}$ can be expressed as the problem of finding a solution to the following equations:
{\setlength\arraycolsep{2pt}
\begin{eqnarray}
\langle \nu_i,\nu_i \rangle &=& 1 \;\;\text{for all}\;\; i=1, \ldots, f, \nonumber \\
\langle \nu_i,\nu_j \rangle &=& -\cos(\pi/n_{ij}) \;\;\text{if faces $F_i$ and $F_j$ are adjacent in $P$}. \label{eqn:hyp}
\end{eqnarray}}
We call these equations (\ref{eqn:hyp}) the \emph{hyperbolic equations}.

To compare this with Vinberg's equations, first note that $P$ is defined by the system of linear inequalities
\begin{displaymath}
 \alpha_i \leq 0, \hspace{3mm} i=1,\ldots, f, \;\;\text{and}\;\; x_1=1,
\end{displaymath}
where the linear functional $\alpha_i \in V^*$ is dual to $\nu_i$ under the Lorentzian inner product.
In other words, $\alpha_i(x)=\langle \nu_i, x \rangle$.

The hyperbolic reflection in the face $F_i$
is given by $$R_i(x) = x - 2 \langle \nu_i, x \rangle \nu_i = x - \alpha_i(x) b_i$$ where the reflection vector is
$b_i = 2 \nu_i$. So the reflection point $[b_i]=[\nu_i]$ in $\R P^n$ is the projective dual of the hyperplane
containing the face $F_i$ with respect to the sphere at infinity in the Klein model of $\H^n$. There is also a well-known geometric construction corresponding to this kind of duality, see for example \cite[p.71]{Thurston_book}.

Thus taking $b_i = 2 \nu_i$ gives a point $t=\{t_i\}=\{2\nu_i\}$ in $\Phi_{\hat P}^{-1}(0)$
corresponding to the hyperbolic structure on $P$.
This follows since if faces $F_i$ and $F_j$ are adjacent in $P$ then
\begin{displaymath}
a_{ij}=\alpha_i(2\nu_j) = 2\langle \nu_i,\nu_j \rangle = -2\cos(\pi/n_{ij}),
\end{displaymath}
and thus
{\setlength\arraycolsep{2pt}
\begin{eqnarray*}
a_{ii}&=&\alpha_i(\nu_i) = 2\langle \nu_i,\nu_i \rangle=2 \;\;\text{for all}\;\; i=1,\ldots,f,\\
a_{ij}a_{ji}&=& 4\cos^2(\pi/n_{ij}) \;\;\text{if $F_i$ and $F_j$ are adjacent in $P$ and $n_{ij}>2$},\\
a_{ij}&=& 0 \:\text{and}\: a_{ji}=0 \;\;\text{if $F_i$ and $F_j$ are adjacent in $P$ and $n_{ij}=2$}.
\end{eqnarray*}}

\begin{proposition}\label{prop:hyperbolicpt} Let $P$ be a hyperbolic Coxeter polyhedron.
Then  the space $\Phi_{\hat P}^{-1}(0)$ contains a single point
corresponding to the hyperbolic structure on $\hat P$.
\end{proposition}
\begin{proof}
The hyperbolic reflection in each face of $P$ is determined by the face and
a reflection point that is dual to the face in the Klein model.
Since $P$ is fixed, the reflection vectors are determined up to scalar multiplication.
By the normalization conditions in $\Phi_{\hat P}$, we see that the reflection vectors
are uniquely determined. Hence the hyperbolic structure on $\hat P$ corresponds to a single point.
\end{proof}

If we solve the hyperbolic equations (\ref{eqn:hyp}) directly, then we obtain many algebraic solutions $\nu_i$.
However when $n \ge 3$, Mostow-Prasad Rigidity shows that there is only one solution
(up to hyperbolic isometries) with geometric meaning.
To find this, we need to check that the $\nu_i$'s give the desired $n$-dimensional convex hyperbolic polytope.

\begin{theorem} {\rm (Vinberg \cite{Vinberg1985})} \label{thm:vinberg2}
Let the Gram matrix $G$ of the set $S$ of vectors $\{\nu_1, \ldots, \nu_f\}$ be an indecomposable matrix (i.e. it cannot be represented as a direct sum of two matrices) with $1$'s along the diagonal and non-positive entries off it. Assume that $S$ spans $V$ and the cone $K$ defined by the inequalities $\langle x, \nu_i \rangle \leq 0 \;(i=1,\ldots,f)$ intersects the Klein model for $\mathbb{H}^n$.
Then $G$ is the Gram matrix of the convex polytope
$P=K \cap \mathbb{H}^n$
bounded by the hyperplanes $H_i=\{x \in \mathbb{H}^n\;|\;\langle x, \nu_i \rangle = 0\}$.
\end{theorem}

The following observation will be useful for computational purposes, when we need to select the
correct geometric solution from the many algebraic solutions to the hyperbolic equations (\ref{eqn:hyp}).
\begin{lemma}\label{lem:choosing}
If the first entry of each $\nu_i$ is non-negative, then the cone $K$ defined by the inequalities
$\langle x, \nu_i \rangle \leq 0 \;(i=1,\ldots,f)$
intersects $\mathbb{H}^n$.
\end{lemma}
\begin{proof}
$x=(1,0,0, \ldots, 0)$ satisfies the inequalities $\langle x, \nu_i \rangle = -\nu_{i1} \leq 0$.
\end{proof}

\begin{remark}\label{rmk:roeder}
In general, it is difficult to find an exact algebraic solution to the hyperbolic equations (\ref{eqn:hyp}). However, in 3-dimensions, Roeder's Matlab program \cite{Roeder2007} can be used to obtain numerical solutions. His construction uses Newton's method and a homotopy to follow the concrete existence proof given by Andreev
(as modified in \cite{RHD}).  Heard's program ``Orb'' \cite{orb} can also be used to
numerically compute hyperbolic structures on the orientable 3-orbifold obtained by doubling a Coxeter polyhedron along its boundary.
In this paper, we will find many exact solutions using Mathematica.
\end{remark}

\begin{remark}
When describing examples in this paper,
we will sometimes abuse notation and identify $V^*$ with $V$ as follows:
If $\alpha_i \in V^*$ has coordinates $(\alpha_{i1},\alpha_{i2}, \ldots ,\alpha_{i,n+1})$,
and $b_j \in V$ has coordinates $(b_{j1},b_{j2}, \ldots ,b_{j,n+1})$ then
$\alpha_i(b_j)= \alpha_{i1} b_{j1} +  \alpha_{i2} b_{j2} +  \ldots +  \alpha_{i,n+1} b_{j,n+1} =\alpha_i \cdot b_j $,
where $\cdot$ denotes the usual Euclidean dot product in $V=\R^{n+1}$.
Faces of a polytope will always be specified by the coordinate vectors
of the corresponding linear functionals.
In particular, for a hyperbolic polytope in $\H^n$, a face with Lorentzian unit normal $\nu_i$
corresponds to the linear functional with
coordinate vector $\alpha_i = J \nu_i$, where
$J \in SL_{\pm}(n+1,\R)$ is the diagonal matrix with diagonal entries $-1,1, \ldots,1$.
\end{remark}

\section{Local restricted deformation spaces of real projective structures near hyperbolic structures}\label{loc_defs}

We now concentrate on the case of a $3$-dimensional Coxeter orbifold $\hat P$.
Recall that real projective structures in the restricted deformation space of
$\hat P$ correspond to solutions to Vinberg's equations (\ref{eqn:normalize})--(\ref{eqn:Vinberg2}).
In \S  \ref{subsec:Z_tang_space} we study the Zariski tangent space to this solution space, and in
\S  \ref{Z_tang_hyperbolic} we study the
Zariski tangent space to the solution space of the hyperbolic equations (\ref{eqn:hyp}).
In \S  \ref{subsec:Thm1_proof} we compare these tangent spaces and combine this with
Garland-Raghunathan-Weil infinitesimal rigidity (\cite{GarlandRag}, \cite{Weil1962})
to prove Theorem \ref{thm:main}. In \S  \ref{subsec:prisms}
we also use these techniques to study the restricted deformation spaces of 3-dimensional compact hyperbolic prisms.

\subsection{The Zariski tangent space to the Vinberg equations}\label{subsec:Z_tang_space}

We now study the Zariski tangent space to the solution space of Vinberg's equations, using the notation from
\S \ref{subsec:proj_rep_space}.
Let $V=\R^4$, let $P$ be a  convex polyhedron in ${\mathbb S}^3$,
and let $\hat P$ be a corresponding Coxeter orbifold.
We assume that $P$ has $f$ faces, and that each
linear functional $\alpha_i \in V^*$ is {\em fixed} for $i=1, \ldots, f$.
Then we have variables $b_i \in V$ for $i=1, \ldots, f$, and the equations have the form
\begin{itemize}
\item $\Phi_{ii} = \alpha_i(b_i) - 2 = 0$,
\item $\Phi_{ij} = \alpha_i(b_j) \alpha_j(b_i) - c_{ij} = 0$ where $c_{ij}$ is a constant
if $n_{ij} \ne 2$,
\item $\Phi^1_{ij} = \alpha_i(b_j)=0$ and $\Phi^2_{ij} = \alpha_j(b_i)=0$ if $n_{ij}=2$.
\end{itemize}

Let $\pi_i : V^f \to V$ denote the projection onto the $i$th factor.
Then the derivative of $\Phi_{ij}$ at $b=(b_1, \ldots, b_f)$, considered as a linear map, is given by:
$$D\Phi_{ij}(\dot b)= \alpha_i(\dot b_j) \alpha_j(b_i)+ \alpha_i(b_j) \alpha_j(\dot b_i)$$
for $\dot b = (\dot b_1, \dot b_2, \ldots, \dot b_f) \in V^f$, or
$$D\Phi_{ij} = a_{ji} (\alpha_i \circ \pi_j) + a_{ij} (\alpha_j \circ \pi_i).$$
Similarly,
$$D\Phi_{ii} = \alpha_i \circ \pi_i, ~D \Phi^1_{ij} = \alpha_i \circ \pi_j \text{ and }  D \Phi^2_{ij} = \alpha_j \circ \pi_i.$$

More explicitly, combining Vinberg's equations gives a function $\Phi_{\hat P} : V^{f} = \R^{4f} \to \R^N$ and the  rows of the $N \times 4f$ Jacobian matrix $D = [D \Phi_{\hat P}]$ are
made up of blocks, each consisting of four entries:
    \begin{align}
  [D\Phi_{ii}]  = & \quad (0, \ldots, 0, \alpha_{i1}, \alpha_{i2}, \alpha_{i3}, \alpha_{i4}, 0, \ldots, 0) \nonumber\\  = &  \quad (0, \ldots, 0, \quad   \underbrace{\alpha_i}_{i \text{th block}} \quad  ,0, \ldots, 0), \nonumber
    \end{align}

     \begin{align}
[D\Phi_{ij}] =  & \quad (0, \ldots, 0, \underbrace{a_{ij}\alpha_j}_{i \text{th block}} ,0, \ldots, 0, \underbrace{a_{ji}\alpha_i}_{j \text{th block}}, 0, \ldots, 0) , \text{ for } n_{ij} \ne 2, \nonumber
    \end{align}

$$[D \Phi^1_{ij}]  = \quad (0, \ldots, 0, \underbrace{0}_{i \text{th block}} ,0, \ldots, 0, \underbrace{\alpha_i}_{j \text{th block}}, 0, \ldots, 0) \text{ for } n_{ij} = 2,$$
and
$$[D \Phi^2_{ij}] = \quad (0, \ldots, 0, \underbrace{\alpha_j}_{i \text{th block}} ,0, \ldots, 0, \underbrace{0}_{j \text{th block}}, 0, \ldots, 0) \text{ for } n_{ij} = 2.$$

Note that this Jacobian matrix has two rows for each edge of $P$ with $n_{ij}=2$, but only one row for each edge with $n_{ij} \geq 3$. (Compare this with \S \ref{subsec:Thm1_proof} below.)

Suppose that $p$ is a point of $\Phi_{\hat P}^{-1}(0)$. Then the  Zariski tangent space at $p$ is the kernel of the Jacobian matrix $D$ evaluated at $p$. We call this the \emph{infinitesimal restricted deformation space of $P$} at $p$ because of
Theorem \ref{isotopy_equiv}.

The next result now follows from the implicit function theorem.

\begin{proposition}\label{prop:fullrank} Let $D = [D\Phi_{\hat P}]$ be the Jacobian matrix of Vinberg's equations for $\hat P$ at $p$. \\
If $4f-N > 0$ and $D$ has full rank, i.e. $\text{rank}\, D = \min(4f,N)$,
then $\mathfrak{D}_P(\hat{P})$
is locally a smooth manifold of dimension $4f-N$ near $p$.
So if $p$ is the hyperbolic point, the hyperbolic structure on the Coxeter 3-orbifold $\hat P$ deforms relative to the mirrors to a real projective structure which is not a hyperbolic structure.\\
If $4f-N \leq 0$ and $D$ has full rank, then $p$ is a isolated point in $\mathfrak{D}_P(\hat{P})$. So if $p$ is the hyperbolic point, the hyperbolic structure on
$\hat P$ is projectively rigid relative to the mirrors in $\mathfrak{D}_P(\hat{P})$.
\end{proposition}
Note that $4f-N=4f-(f+e+e_2)=3f-e-e_2$.
The results in  \cite{Choi2006} are obtained by showing that $D$ has full rank in
the \emph{orderable} case.

The following example illustrates
 the role of orderability; in this case,
permutations of rows of $D$ are sufficient to show that $D$ has full rank.
This example was originally studied by Benoist \cite{Benoist2006}  and is orderable.

Here and throughout the paper, we use the following notation. Given a diagram of
a 3-dimensional hyperbolic polyhedron, if an edge is labelled $e_i$, then its dihedral angle is $\pi/e_i$. Moreover, $\alpha_i$ is the linear functional defining the face $F_i$.

\begin{example}\label{exp:prism}
Figure 1 shows a 3-dimensional compact triangular prism $P$. This satisfies the conditions (A1)--(A4) of Andreev's Theorem so defines a hyperbolic Coxeter orbifold $\hat P$.

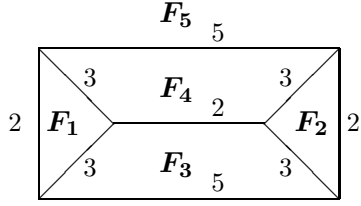
\begin{figure}[h]
\centering
\setlength{\unitlength}{1cm}
\begin{picture}(4,2.7)
\put(0.1,0.9){\large $\boldsymbol{F_1}$}
\put(3.4,0.9){\large $\boldsymbol{F_2}$}
\put(1.6,2.35){\large $\boldsymbol{F_5}$}
\put(1.6,1.35){\large $\boldsymbol{F_4}$}
\put(1.6,0.35){\large $\boldsymbol{F_3}$}
\put(-0.4,0.9){$2$}
\put(0,0){\line(0,1){2}}
\put(0.6,0.3){$3$}
\put(0,0){\line(1,1){1}}
\put(0.6,1.5){$3$}
\put(0,2){\line(1,-1){1}}
\put(2.3,2.1){$5$}
\put(0,2){\line(1,0){4}}
\put(2.3,1.1){$2$}
\put(0,0){\line(1,0){4}}
\put(2.3,0.1){$5$}
\put(4,0){\line(0,1){2}}
\put(3.2,1.5){$3$}
\put(4,0){\line(-1,1){1}}
\put(3.2,0.3){$3$}
\put(4,2){\line(-1,-1){1}}
\put(4.1,0.9){$2$}
\put(1,1){\line(1,0){2}}
\end{picture}
\caption{A compact hyperbolic triangular prism}
\end{figure}

The following table shows that $\hat P$ is orderable:
\begin{displaymath}
\begin{array}{|c||c|c|}
\hline
  & \text{faces} & \text{faces of higher index} \\
  & \text{sharing edges of order 2} & \text{not sharing edges of order 2}\\
\hline
F_1 & F_5 & F_3, \; F_4 \\
\hline
F_2 & F_5 & F_3, \; F_4 \\
\hline
F_3 & F_4 & F_5 \\
\hline
F_4 & F_3 & F_5 \\
\hline
F_5 & F_1,\; F_2 & \varnothing \\
\hline
\end{array}
\end{displaymath}
The $17 \times 20$ Jacobian matrix $D=[D\Phi_{\hat P}]$ is shown on the left below. By permuting the rows,
we obtain the matrix on the right. (The new ordering corresponds to the entries in the rows of the table above.)
{\tiny
\begin{displaymath}
\
D =
\left( \begin{array}{ccccc}
\alpha_{1} & 0 & 0 & 0 & 0 \\
0 & \alpha_{2} & 0 & 0 & 0 \\
0 & 0 & \alpha_{3} & 0 & 0 \\
0 & 0 & 0 & \alpha_{4} & 0 \\
0 & 0 & 0 & 0 & \alpha_{5} \\
a_{13}\alpha_{3} & 0 & a_{31}\alpha_{1} & 0 & 0 \\
a_{14}\alpha_{4} & 0 & 0 & a_{41}\alpha_{1} & 0 \\
\alpha_{5} & 0 & 0 & 0 & 0 \\
0 & 0 & 0 & 0 & \alpha_{1} \\
0 & a_{23}\alpha_{3} & a_{32}\alpha_{2} & 0 & 0 \\
0 & a_{24}\alpha_{4} & 0 & a_{42}\alpha_{2} & 0 \\
0 & \alpha_{5} & 0 & 0 & 0 \\
0 & 0 & 0 & 0 & \alpha_{2} \\
0 & 0 & \alpha_{4} & 0 & 0 \\
0 & 0 & 0 & \alpha_{3} & 0 \\
0 & 0 & a_{35}\alpha_{5} & 0 & a_{53}\alpha_{3} \\
0 & 0 & 0 & a_{45}\alpha_{5} & a_{54}\alpha_{4}
\end{array} \right)
\sim
\left( \begin{array}{ccccc}
\boldsymbol{\alpha_{1}} & 0 & 0 & 0 & 0 \\
\boldsymbol{\alpha_{5}} & 0 & 0 & 0 & 0 \\
\boldsymbol{a_{13}\alpha_{3}} & 0 & a_{31}\alpha_{1} & 0 & 0 \\
\boldsymbol{a_{14}\alpha_{4}} & 0 & 0 & a_{41}\alpha_{1} & 0 \\
0 & \boldsymbol{\alpha_{2}} & 0 & 0 & 0 \\
0 & \boldsymbol{\alpha_{5}} & 0 & 0 & 0 \\
0 & \boldsymbol{a_{23}\alpha_{3}} & a_{32}\alpha_{2} & 0 & 0 \\
0 & \boldsymbol{a_{24}\alpha_{4}} & 0 & a_{42}\alpha_{2} & 0 \\
0 & 0 & \boldsymbol{\alpha_{3}} & 0 & 0 \\
0 & 0 & \boldsymbol{\alpha_{4}} & 0 & 0 \\
0 & 0 & \boldsymbol{a_{35}\alpha_{5}} & 0 & a_{53}\alpha_{3} \\
0 & 0 & 0 & \boldsymbol{\alpha_{4}} & 0 \\
0 & 0 & 0 & \boldsymbol{\alpha_{3}} & 0 \\
0 & 0 & 0 & \boldsymbol{a_{45}\alpha_{5}} & a_{54}\alpha_{4} \\
0 & 0 & 0 & 0 & \boldsymbol{\alpha_{5}} \\
0 & 0 & 0 & 0 & \boldsymbol{\alpha_{1}} \\
0 & 0 & 0 & 0 & \boldsymbol{\alpha_{2}}
\end{array} \right)
\end{displaymath}}

Note that every coefficient $a_{ij}$ appearing in the Jacobian matrix $D$ corresponds to an edge with
$n_{ij}\geq  3$, so satisfies $a_{ij}<0$. Thus Lemma \ref{lem:ind} below
implies that the following submatrices 
have full rank.
(Note that the orderability of $\hat P$ is used here.)

{
\begin{displaymath}
\left( \begin{array}{c}
\boldsymbol{\alpha_1} \\
\boldsymbol{\alpha_5} \\
\boldsymbol{a_{13}\alpha_3} \\
\boldsymbol{a_{14}\alpha_4}
\end{array} \right),
\left( \begin{array}{c}
\boldsymbol{\alpha_2} \\
\boldsymbol{\alpha_5} \\
\boldsymbol{a_{23}\alpha_3} \\
\boldsymbol{a_{24}\alpha_4}
\end{array} \right),
\left( \begin{array}{c}
\boldsymbol{\alpha_3} \\
\boldsymbol{\alpha_4} \\
\boldsymbol{a_{35}\alpha_5}
\end{array} \right),
\left( \begin{array}{c}
\boldsymbol{\alpha_4} \\
\boldsymbol{\alpha_3} \\
\boldsymbol{a_{45}\alpha_5}
\end{array} \right),
\left( \begin{array}{c}
\boldsymbol{\alpha_5} \\
\boldsymbol{\alpha_1} \\
\boldsymbol{\alpha_2}
\end{array} \right).
\end{displaymath}}
Therefore, the Jacobian matrix $D$ has full rank $= 17$, and thus $\mathfrak{D}_P(\hat{P})$ is a smooth manifold of dimension $3$, since the dimension of the null space of $D$ is always $3=20-17$.
\end{example}

The following observation will be used again in \S \ref{subsec:lin_test}.
\begin{lemma}\label{lem:ind}
Let $P$ be a 3-dimensional convex polyhedron in $\mathbb{S}^3$ defined by linear inequalities $\alpha_i \le 0$ where $\alpha_i \in V^*$, and let $F_i$ be the face of $P$ determined by $\alpha_i$. Suppose that the faces
$F_{i_1}$, $F_{i_2}$, $F_{i_3}$ are adjacent to the face $F_{i_4}$.
Then the four linear functionals $\alpha_{i_1}$, $\alpha_{i_2}$, $\alpha_{i_3}$, $\alpha_{i_4}$ are linearly independent.
\end{lemma}
\begin{proof}
If  not, then the linear functionals $\alpha_{i_1}$, $\alpha_{i_2}$, $\alpha_{i_3}$, $\alpha_{i_4}$ lie in a codimension one subspace of $V^*$ determined  by a non-zero vector $b \in V$. In other words,
\begin{displaymath}
\alpha_{i_1}(b)=\alpha_{i_2}(b)=\alpha_{i_3}(b)=\alpha_{i_4}(b)=0.
\end{displaymath}
Furthermore, we know that three linear inequalities and one equality
\begin{displaymath}
\alpha_{i_1} \leq 0,\; \alpha_{i_2} \leq 0,\; \alpha_{i_3} \leq 0,\; \alpha_{i_4} = 0
\end{displaymath}
give a (2-dimensional) triangle since the three faces
$F_{i_1}$, $F_{i_2}$, $F_{i_3}$  meet the plane containing the face $F_{i_4}$ in lines. Moreover, these three lines have no common intersection point. This is a contradiction since the three lines meet at the point $b$ by our assumption.
\end{proof}

\subsection{The Zariski tangent space to the hyperbolic equations}\label{Z_tang_hyperbolic}

 Assume $P$ is a finite volume 3-dimensional hyperbolic polyhedron where the dihedral angle at an edge $e_{ij}$
 equals $\pi/{n_{ij}}$ for an integer $n_{ij} \geq 2$. Andreev's theorem  \cite{Andreev19702} characterizes such polyhedra.

Constructing such a hyperbolic polyhedron $P$
is the same as solving the system
of hyperbolic equations (\ref{eqn:hyp}) for the unit normals $\nu_i$ from \S \ref{subsec:hyp_point}.
Equivalently we can write these equations
in terms of the reflection vectors $b_i = 2\nu_i$.  This gives
the following system of $n=f+e$  equations:
\begin{displaymath}
 \Psi_{ii} =  \langle b_i, b_i \rangle - 4 = 0 \;\; \text{and} \;\;
 \Psi_{ij} = \langle b_i,\ b_j \rangle + 4 \cos(\pi/n_{ij}) = 0.
 \end{displaymath}
 Combining these gives a function
  $\Psi_{\hat P}: V^f = \R^{4f} \to \R^n$ and $\Psi_{\hat P}^{-1}(0)$
contains the convex polyhedra in $\H^3$ with the desired dihedral angles.
By Andreev's Theorem (or Mostow-Prasad rigidity), there is a unique such polyhedron up to hyperbolic isometries;
this corresponds to a 6-dimensional manifold contained in $\Psi_{\hat P}^{-1}(0)$.

Now consider the
derivative $D\Psi_{\hat P}$ at a hyperbolic point $t$.
If $\alpha_i = \langle \nu_i, \cdot \rangle$ are the linear functionals defining the faces
of the hyperbolic polyhedron then
$$D \Psi_{ij} (\dot b) = \langle \dot b_i, b_j \rangle + \langle  b_i,  \dot b_j \rangle
 = 2 \alpha_j (\dot b_i) + 2 \alpha_i(\dot b_j),$$
 or $$D\Psi_{ij} = 2 \alpha_j \circ \pi_i + 2 \alpha_i \circ \pi_j.$$
  When $i=j$ this becomes
 $$D\Psi_{ii}= 4 \alpha_i \circ \pi_i.$$
Equivalently, the rows of the $n \times 4f$ Jacobian matrix $\hat{D} = [D \Psi_{\hat P}]$ are
made up of blocks, each consisting of four entries:
    \begin{align}
    [D\Psi_{ii}]  = & \quad (0, \ldots, 0, 4\alpha_{i1}, 4\alpha_{i2}, 4\alpha_{i3}, 4\alpha_{i4}, 0, \ldots, 0) \nonumber\\  = &  \quad (0, \ldots, 0, \quad \quad \: \underbrace{4\alpha_i}_{i \text{th block}} \quad \quad \: ,0, \ldots, 0) \nonumber,
    \end{align}
and
    \begin{align}
   [D\Psi_{ij}]  = & \quad (0, \ldots, 0, \underbrace{2\alpha_j}_{i \text{th block}} ,0, \ldots, 0, \underbrace{2\alpha_i}_{j \text{th block}}, 0, \ldots, 0) \nonumber.
    \end{align}
Then the Zariski tangent space to $\Psi_{\hat P}^{-1}(0)$ at $t$ is $\ker D\Psi_{\hat P}$.

\subsection{The proof of Theorem \ref{thm:main}}\label{subsec:Thm1_proof}

We now assume that $P$ is a convex ideal polyhedron in $\H^3$ with all edges of order $3$. Then all vertices are trivalent, and we have assumed that $P$ is not a tetrahedron.  So it follows from Lemma \ref{lem:kp} that $k(P)=0$,
and that the results from \S \ref{subsec:def_vs_hol} apply.

To prove Theorem 1,
we use the results from \S \ref{subsec:Z_tang_space} and \S \ref{Z_tang_hyperbolic}
to compare the Jacobian matrices $D=[D\Phi_{\hat P}]$ for real projective structures
and $\hat D = [D\Psi_{\hat P}]$ for hyperbolic structures.
Since $P$ contains no edges of order 2 we have $N=n$, and each  $a_{ij}$ is non-zero.
Further,  $a_{ij}=a_{ji}$ at a hyperbolic point in $\mathfrak{D}_P(\hat{P})$.
Hence, each row of $D$ is a non-zero scalar multiple of a row of $\hat{D}$,
so the ranks of $D$ and $\hat{D}$ are equal.

We now use the infinitesimal rigidity of the hyperbolic structure on $\hat P$ to compute the rank of $\hat D$.
The arguments from Proposition \ref{prop:eqnsol} show that $\Psi_{\hat P}^{-1}(0)$
is locally isomorphic to the algebraic variety $\text{Hom} (\pi_1(\hat P), O_0(1,3))$
near a hyperbolic point with holonomy representation $h_0$.
Thus, by the work of Weil  \cite{Weil1964}, the Zariski tangent space to $\Psi_{\hat P}^{-1}(0)$ at this point
corresponds to the space of 1-cocyles in group cohomology
$Z^1(\pi_1(\hat P), so(1,3)_{Ad})$, where $\pi_1(\hat P)$ acts on $so(1,3)$
via the representation $Ad \circ h_0$. (See also \cite{Raghunathan}, \cite{Kapovich}, \cite{Kapovich1994}.)

Let $Q$ be the compact orbifold obtained by truncating the cusps of $\hat P$; then $Q$ has a boundary $\bd Q$ consisting of $(3,3,3)$-triangle orbifolds.
Now consider the exact sequence
$$H^1(Q,\bd Q; so(1,3)_{Ad}) \to H^1(Q; so(1,3)_{Ad})  \to H^1(\bd Q; so(1,3)_{Ad}).$$
By Garland-Raghunathan-Weil infinitesimal rigidity (\cite{GarlandRag}, \cite{Weil1962})
the parabolic group cohomology $PH^1(\pi_1(\hat P); so(1,3)_{Ad})=0$.
This implies that
$$H^1(Q,\bd Q; so(1,3)_{Ad})=0.$$ Further, the $(3,3,3)$-triangle group is infinitesimally rigid in $O_0(1,3)$
so $$H^1(\bd Q; so(1,3)_{Ad})=0.$$ Hence  $$H^1(\pi_1(\hat P); so(1,3)_{Ad}) \cong H^1(Q; so(1,3)_{Ad}) =0.$$
It follows that $$\dim \ker D\Psi_{\hat P} = \dim Z^1(\pi_1(\hat P), so(3,1)_{Ad}) = \dim so(1,3) =6.$$
In fact, Weil's argument in \cite{Weil1964} shows that a neighbourhood of $h_0$ in
 $\text{Hom}(\pi_1(\hat P), O_0(1,3))$
 coincides with the orbit of $h_0$ under the group of hyperbolic isometries, and
 that this is a  locally a smooth $6$-manifold since the hyperbolic holonomy group
 $h_0(\pi_1(\hat P))$
has trivial centralizer.

Since all vertices of $P$ are trivalent  we have $3v=2e$, and  since $v-e+f=2$ it follows that $4f-N=3f-e= 6 = 4f - \rank D\Psi_{\hat P}$.
Hence $\rank D\Phi_{\hat P} = \rank D\Psi_{\hat P} = N$ and $D\Phi_{\hat P}$ has full rank.
 Therefore, by Proposition \ref{prop:fullrank}, a neighbourhood of $t$ in $\mathfrak{D}_P(\hat{P})$ is a smooth $6$-dimensional manifold.
 This completes the proof of Theorem \ref{thm:main}.

\begin{remark} \label{rmk:solutions}
Looking more closely at the proof of Theorem \ref{thm:main},
we see that $\ker D\Phi_{\hat P}$ is \emph{equal} to  $\ker D \Psi_{\hat P}$ and is given by the tangent space to orbit of
$b=(b_1, \ldots,b_f)$ under the $6$-dimensional group of hyperbolic isometries.
This gives a very nice geometric interpretation of the infinitesimal deformations relative to the mirrors
provided by Theorem \ref{thm:main}. Consider a convex hyperbolic polyhedron $P$ in the
Klein model with $P \subset \H^3  \subset \R P^3$.
The reflection in each face $F_i$ has a fixed point $[b_i]$ outside the sphere at infinity,
corresponding to the reflection vector $b_i$ for the face. Let $g_t$ be a 1-parameter family of
isometries of $\H^3$ with $g_0 = identity$. Then
$b(t) =  (g_t(b_1), \ldots, g_t(b_f))$ is a curve in $V^f=\R^{4f}$ whose
derivative at the identity $\dot b =  (\dot g(b_1), \ldots, \dot g(b_f))$
is in the kernel of $D\Phi_{\hat P}$.  Here $\dot g \in so(1,3)$ is an infinitesimal isometry of $\H^3$.
In other words, all the infinitesimal projective deformations relative to the mirrors
are obtained by fixing the polyhedron faces and moving the fixed points of face reflections by
infinitesimal hyperbolic isometries. It would be very interesting to extend this observation to give an explicit description
of the local projective deformations relative to the mirrors.
\end{remark}

\begin{remark}\label{rem:higherorder}
This argument extends to convex hyperbolic polyhedra $P$ with trivalent but possibly hyperinfinite vertices, provided {\em all edges have order at least $3$} and $k(P)=0$.
In general, such a polyhedron is non-compact of infinite volume, but can be truncated along planes orthogonal to the faces at each hyperinfinite vertex to give a compact convex polyhedron. Again infinitesimal rigidity applies since all vertex cross sections give hyperbolic triangle groups, hence are rigid.
The argument given above then shows that the restricted deformation space $\mathfrak{D}_P(\hat{P})$ is again locally a smooth $6$-dimensional manifold, provided
{\em all edges have order at least $3$}. So ``almost all'' assignments of orders to the edges of $P$ give a hyperbolic Coxeter orbifold $\hat P$
that is projectively deformable relative to the mirrors.
\end{remark}

\subsection{Deformations of prisms}\label{subsec:prisms}

Assume that $P$ is a $3$-dimensional compact hyperbolic Coxeter polyhedron, and $\hat P$ a corresponding Coxeter orbifold.
Theorem \ref{thm:orderable7} below shows that whenever $\hat{P}$ is orderable and the number of faces of $P$ is greater than $7$,
$\hat P$ is projectively rigid relative to the mirrors. Thus, there are only finitely many combinatorial types of convex hyperbolic
polyhedra with orderable compact Coxeter 3-orbifold structures that are
projectively deformable relative to the mirrors. However when $\hat{P}$ is not orderable,  Proposition \ref{prop:prism} shows that this is no longer true.

\begin{theorem}\label{thm:orderable7}
Let $P$ be a $3$-dimensional compact hyperbolic Coxeter polyhedron. Suppose that $\hat{P}$ is orderable. If the number of faces of $P$ is greater than $7$, then it is projectively rigid relative to the mirrors in $\mathfrak{D}_P(\hat{P)}$.
\end{theorem}
\begin{proof} Since $P$ is compact, every vertex is trivalent and is
adjacent to an edge of order 2. (See, for example, Andreev's condition (A1) in \S \ref{subsec:andreev}.)
Also, $v > 10$, since $3v=2e$ and $v-e+f=2$. This implies that $e_2>5$, and thus $3f-e-e_2=6-e_2 \leq 0$.
Hence, there is no local deformations of $P$ relative to the mirrors by Theorem \ref{thm:choi}
\end{proof}

\begin{proposition}\label{prop:prism}
For any natural number $f \geq 7$, there exists a $3$-dimensional compact hyperbolic prism $P$ with $f$ faces and a corresponding Coxeter orbifold $\hat{P}$ that can be projectively deformed relative to the mirrors.
\end{proposition}
\begin{proof}
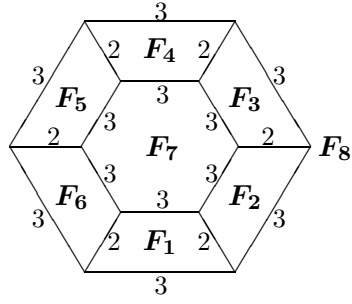
\begin{figure}[h]
\centering
\setlength{\unitlength}{10mm}
\begin{picture}(5,4)
\put(2.27,0.5){\large $\boldsymbol{F_1}$}
\put(3.4,1.1){\large $\boldsymbol{F_2}$}
\put(3.4,2.4){\large $\boldsymbol{F_3}$}
\put(2.27,3.05){\large $\boldsymbol{F_4}$}
\put(1.1,2.4){\large $\boldsymbol{F_5}$}
\put(1.1,1.1){\large $\boldsymbol{F_6}$}
\put(2.3,1.75){\large $\boldsymbol{F_7}$}
\put(4.6,1.75){\large $\boldsymbol{F_8}$}
\put(2.45,2.45){$3$}
\put(2.45,1.05){$3$}
\put(3.1,1.4){$3$}
\put(3.1,2.1){$3$}
\put(1.75,1.4){$3$}
\put(1.75,2.1){$3$}
\put(3.85,1.9){$2$}
\put(1.0,1.9){$2$}
\put(1.8,3.05){$2$}
\put(3.0,3.05){$2$}
\put(1.8,0.5){$2$}
\put(3.0,0.5){$2$}
\put(2.43,3.56){$3$}
\put(2.43,-0.09){$3$}
\put(4.0,2.7){$3$}
\put(4.0,0.8){$3$}
\put(0.8,2.7){$3$}
\put(0.8,0.8){$3$}

\put(1.98,2.74){\line(1,0){1.03}}
\put(3.02,2.74){\line(3,-5){0.52}}
\put(1.98,2.74){\line(-3,-5){0.52}}
\put(1.98,1){\line(1,0){1.03}}
\put(3.02,1){\line(3,5){0.52}}
\put(1.98,1){\line(-3,5){0.52}}

\put(1.5,3.54){\line(1,0){2}}
\put(3.5,3.54){\line(3,-5){1.0}}
\put(1.5,3.54){\line(-3,-5){1.0}}
\put(1.5,0.2){\line(1,0){2}}
\put(3.5,0.2){\line(3,5){1.0}}
\put(1.5,0.2){\line(-3,5){1.0}}

\put(1.5,0.2){\line(3,5){0.48}}
\put(3.5,0.2){\line(-3,5){0.48}}
\put(4.5,1.87){\line(-1,0){0.96}}
\put(0.5,1.87){\line(1,0){0.96}}
\put(1.5,3.54){\line(3,-5){0.48}}
\put(3.5,3.54){\line(-3,-5){0.48}}
\end{picture}
\caption{A hexagonal prism with a rotational symmetry about the axis passing through the centers of faces $F_7$ and $F_8$}
\end{figure}
Let $n$ be a natural number greater than $4$.
We will construct a prism $P$ with $f = n+2$ faces such that edge orders of the top and bottom $n$-gons of $P$ are $3$ and all the remaining edges orders are $2$.
(Figure 2 shows the case where $n=6$.)
Note that $P$ is not orderable.

Suppose that $c=\cos(2\pi/n)$ and $s=\sin(2\pi/n)$. Using the rotational symmetry of $P$, we can find a compact $3$-dimensional hyperbolic prism as follows. Let the unit normals $\nu_i$ have coordinate vectors:

{\setlength\arraycolsep{2pt}
\begin{eqnarray*}
\nu_1 & = & \big(-\sqrt{ \frac{c}{1-c}}, \sqrt{ \frac{1}{1-c}},\:0,\:0 \big),\\
\nu_{k+1}&=&L^{k}\nu_1 \:\:\text{with}\:\: k=1, \ldots, n-1,\\
\nu_{n+1} &=& \big(-\sqrt{ \frac{1-c}{4c}},\:0,\:0, \sqrt{ \frac{1+3c}{4c}}\big),\\
\nu_{n+2}&=&\big(-\sqrt{ \frac{1-c}{4c}},\:0,\:0, -\sqrt{ \frac{1+3c}{4c}}\big),
\end{eqnarray*}}
where
\begin{displaymath}
L = \left( \begin{array}{cccc}
1 & 0 & 0 & 0\\
0 & c & -s & 0 \\
0 & s & c & 0 \\
0 & 0 & 0 & 1
\end{array} \right).
\end{displaymath}

Then
{\setlength\arraycolsep{2pt}
\begin{eqnarray*}
& & \langle \nu_i, \nu_i \rangle = 1 \:\:\text{for all}\:\: i=1, \ldots, n+2,\\
& & \langle \nu_{n+1}, \nu_k \rangle = -1/2 \:\:\text{and}\:\: \langle \nu_{n+2}, \nu_k \rangle = -1/2  \:\:\text{for all}\:\: k=1, \ldots, n,\\
& & \langle \nu_1, \nu_{2} \rangle = \langle \nu_2, \nu_{3} \rangle = \cdots = \langle \nu_{n-1}, \nu_{n} \rangle = \langle \nu_{n}, \nu_{1} \rangle =0,
\end{eqnarray*}}
where $\langle \cdot\,,\cdot \rangle$ denotes the Lorentzian inner product.

Our aim is to solve Vinberg's equations (\ref{eqn:normalize})--(\ref{eqn:Vinberg2}).
Let the linear functionals $\alpha_i$ defining the faces have coordinate vectors $J\nu_i$.
We restrict ourselves to  solutions of the special form:
{\setlength\arraycolsep{2pt}
\begin{eqnarray*}
b_1 &=& \big(\sqrt{\frac{4c}{1-c}},\: \sqrt{ \frac{4}{1-c}},\:0,\:b_{1,4} \big),\\
b_{k+1} &=& L^{k}b_1 \:\:\text{with}\:\: k=1, \ldots, n-1, \\
b_{n+1} &=& (b_{n+1,1},0,0,b_{n+1,4}) \:\:\text{with}\:\: \alpha_{n+1}(b_{n+1})=2,\\
b_{n+2} &=& (b_{n+2,1},0,0,b_{n+2,4}) \:\:\text{with}\:\: \alpha_{n+2}(b_{n+2})=2.
\end{eqnarray*}}
where $b_{n+1,2}=b_{n+1,3}=b_{n+2,3}=b_{n+2,3}=0$, and we
take $\{ b_{1,4}, b_{n+1,4}, b_{n+2,4}\}$ as \emph{three} free variables.

We now solve Vinberg's equations.
By computation, we obtain
$$a_{kk}=\alpha_k(b_k)=\alpha_1(b_1)=2$$ for all $k=1, \ldots, n$, since $(L^k)^T L^k=I$, where $M^T$ denotes the transpose of the matrix $M$. Then we obtain
{\setlength\arraycolsep{2pt}
\begin{eqnarray*}
a_{j,j+1} &=& \alpha_j(b_{j+1})=\alpha_1(b_2)=-2c/(1-c)+2c/(1-c)=0,\\
a_{j+1,j} &=& \alpha_{j+1}(b_{j})=\alpha_2(b_1)=-2c/(1-c)+2c/(1-c)=0,
\end{eqnarray*}}
for all $j=1, \ldots, n-1$. Similarly, $a_{1,n}=a_{n,1}=0$.
Moreover, if the \emph{two} equations
\begin{equation} \label{2_prism_eqns}
a_{n+1,1} a_{1,n+1}=1 \;\; \text{and} \;\; a_{n+2,1} a_{1,n+2}=1
\end{equation}
are satisfied, then the rotational symmetry implies that {\em all} of Vinberg's equations (\ref{eqn:normalize})--(\ref{eqn:Vinberg2})  are satisfied, since
\begin{displaymath}
a_{n+1,k}=\alpha_{n+1}(b_k)=\alpha_{n+1}(b_1) \;\; \text{and}\;\; a_{k,n+1}=\alpha_{k}(b_{n+1})=\alpha_{1}(b_{n+1}),
\end{displaymath}
for all $k=1, \ldots, n$. Similarly, $a_{n+2,k}=a_{n+2,1}$ and $a_{k,n+2}=a_{1,n+2}$ for all $k=1, \ldots, n$.

Since we have three free variables subject to the two equations (\ref{2_prism_eqns}),
it follows that the hyperbolic Coxeter orbifold $\hat{P}$ has a one parameter family of real projective structures.
(Of course, the dimension of the restricted deformation space might be greater than $1$.)
This completes the proof.
\end{proof}

\section{The numerical and algebraic computations of restricted deformation spaces}\label{computations}

The restricted deformation space $\mathfrak{D}_P(\hat{P)} \cong \Phi^{-1}(0)$ of a Coxeter orbifold $\hat P$
is defined by Vinberg's system of polynomials $\Phi=\Phi_{\hat P}$ and each of these  has total degree $\leq 2$.
It is difficult to understand the general properties of these algebraic varieties. Thus we examine the infinitesimal and local restricted deformation spaces
for some interesting examples of Coxeter 3-orbifolds arising from compact hyperbolic cubes and dodecahedra.
This work uses a combination of several theoretical and computational methods.

In \S  \ref{subsec:main_alg} we outline the main algorithm used for our computations. In \S  \ref{subsec:lin_test} we give a simple linear test for projective rigidity rel mirrors. In \S  \ref{subsec:notation} we describe the notation used in our figures and tables. Finally, in \S \ref{subsec:cubes} and \S  \ref{subsec:dodec}, we provide
detailed tables describing our results on the restricted projective deformation spaces for cubes and dodecahedra,
and give detailed descriptions of the methods used.

\subsection{The main algorithm for computing local restricted deformations}\label{subsec:main_alg}
We use the following steps to compute the local restricted deformation spaces of 3-dimensional compact hyperbolic Coxeter cubes and dodecahedra.
\begin{enumerate}
\item
We tabulate the $3$-dimensional compact hyperbolic cubes (or dodecahedra) satisfying the conditions of Andreev's theorem (A1)-(A4).
(See cu.m and do.m in \cite{Lee}.)
To obtain manageable finite lists, we restrict the possible edge orders as specified in Theorems 2 and 3. This gives us 34 Coxeter orbifolds (\textsf{cu1}-\textsf{cu34})  based on the cube, and 13 Coxeter orbifolds (\textsf{do1}-\textsf{do13}) based on the dodecahedron.

\item We apply the linear test of rigidity in \textsection \ref{subsec:lin_test} by hand.
If the test shows rigidity relative to the mirrors, we stop here and conclude that our orbifold is projectively rigid relative to the mirrors.

\item Next, we explictly construct
the 3-dimensional compact hyperbolic Coxeter cubes (or dodecahedra) obtained in step 1.
To do this we first choose three faces meeting at a vertex and put the normals to these faces
into a standard position. We then
use Mathematica to solve the hyperbolic equations (\ref{eqn:hyp}) for the remaining unit normals.
This gives us explicit linear functionals  $\alpha_i$ defining the hyperbolic polyhedron.

For cubes, it is not difficult to find exact \emph{algebraic} values for $\alpha_i$, since the number of $\alpha_i$ is not large.
However, for dodecahedra it is difficult to find these algebraic values.
By utilizing a rotational symmetry of \textsf{do13}, we find exact algebraic values of $\alpha_i$ for \textsf{do13} by hand.  We then obtain \emph{numerical} values of $\alpha_i$ for the remaining dodecahedra \textsf{do1-do12} by deforming the dihedral angles of \textsf{do13}.
To obtain the numerical values of $\alpha_i$, we utilize Mathematica
where we can adjust the accuracy to make the errors as small as desired.
In fact, we maintain 150 digits of precision in internal computations.

To check that the $\alpha_i$ obtained here give the desired 3-dimensional compact hyperbolic polyhedron, we use Theorem \ref{thm:vinberg2}. In fact, the first coordinate of each $\alpha_i$
obtained by Mathematica is non-negative, so
 we can easily apply Lemma \ref{lem:choosing}

\item We compute the dimension of the Zariski tangent space to $\Phi^{-1}(0)$
 for the hyperbolic point, i.e. the dimension of the null space of the Jacobian matrix
 $D = [D\Phi]$ at the corresponding point. If $D$ is of full rank, step 5 is followed. Otherwise, step 6 is followed.

For cubes, we use the exact algebraic values obtained in step 3.
However, for dodecahedra,   we have only numerical values of $\alpha_i$ for the dodecahedra other than \textsf{do13}.
To see the accuracy of the numerical calculation of $\rank D$
we use the singular value decomposition of the Jacobian matrix $D$.  In general, the best method for determining the rank of a matrix is to count the number of singular values greater than a judiciously chosen positive number
\cite{Trefethen}. Note that the singular values of the matrix are non-negative real numbers. We check the minimum of the singular values of $D$ to determine whether $D$ is of full rank or not.

\item  If $D$ is of full rank, then the dimension of a neighbourhood of the hyperbolic point is determined by the kernel of $D$. That is,
by Proposition \ref{prop:fullrank}, the dimension of the space of infinitesimal restricted deformations is the same as the dimension of the space of local restricted deformations. Therefore, in this case, the algorithm stops and we obtain answers.

\item If $D$ is rank-deficient,
we attempt to obtain the Gr\"{o}bner basis of the ideal $\mathcal I$ generated by $\{\Phi_k=0\}_{k=1}^{N}$
with respect to a lexicographic order on the variables.
First, we choose new coordinates on $\mathbb{R}^{4f}$ by letting $c_i=b_i-t_i$, where $t=\{t_i\}$
corresponds to our hyperbolic point in $\mathfrak{D}_P(\hat{P})$.
In this coordinate system, the hyperbolic point corresponds to the origin $0$.

In general, the entries of $\alpha_i$ are complicated, and thus sometimes it is difficult
to calculate the Gr\"{o}bner basis of the ideal  $\mathcal I$
directly by using Mathematica.
Therefore, we express the entries of $\alpha_i$ as elements in a field $\mathbb{Q}(\theta)$
generated by an algebraic number $\theta$,
to improve the speed of calculations.
In general, the arithmetic within a \emph{fixed} finite extension of $\mathbb{Q}$ is much faster than arithmetic within the field of complex numbers.
Using this technique, we obtain a list of the dimensions of local restricted deformation spaces
for all the cubes.
For \textsf{do13}, using the rotational symmetry, we compute the dimension of its local  restricted deformation space.
Consequently, we get a list of the dimensions of local restricted deformation spaces  for all the dodecahedra.
\end{enumerate}

\subsection{A linear test for rigidity}\label{subsec:lin_test}
The following \emph{linear test for rigidity} provides a simple, direct
proof that seventeen cubes (\textsf{cu1-cu14, cu16, cu20, cu23}) are projectively rigid relative to the mirrors.
For other cubes, we go to the next step of the algorithm.
Hence, the Mathematica files at \cite{Lee} contain only the results for the remaining seventeen cubes.

Let $P$ be a $3$-dimensional Coxeter polyhedron in $\mathbb{S}^3$.
Then there is a simple method to show the rigidity of the corresponding orbifold $\hat P$.
\begin{enumerate}
\item Find all the faces having more than two edges of order 2. We call them the \emph{rigid faces at level 1}.

\item Relabel all edges of rigid faces at level 1 to become edges of order 2.

\item Again, find all other faces having more than two edges of order 2. We call them  \emph{rigid faces at level 2}. Relabel all edges of these faces to become edges of order 2.

\item Continue the process this manner.

\item If every face of $P$ occurs as a rigid face at level $k$ for some $k \ge 1$,
then we can conclude $\hat P$ is projectively rigid relative to the mirrors.
\end{enumerate}
This test is derived from the following two facts.
\begin{itemize}
\item If a face $F_i$ has more than two edges that are of order $2$, say
 $\{F_{ij_1}, \ldots, F_{ij_m}\}$ with $m \geq 3$, then $b_i$
 can be eliminated as a variable. This follows since
 $b_i$ satisfies a system of linear equations
\begin{displaymath}
    \alpha_{j_1}(b_i)=c_{j_1}, \;\;\ldots \;\;, \alpha_{j_m}(b_i)=c_{j_m}
\end{displaymath}
where $c_{j_1},\ldots, c_{j_m}$ are constants, and these
determine $b_i$ uniquely by Lemma \ref{lem:ind}

\item If $b_i$ is no longer a variable, then $\alpha_i(b_j)\alpha_j(b_i)=4\cos^2(\pi/n_{ij})$ is a \emph{linear} equation for $b_j$.
\end{itemize}

As the number of edges of order $2$ increases, this test becomes more effective. In particular, it is often useful if $P$ is a compact $3$-dimensional hyperbolic Coxeter polyhedron.

\begin{remark}
Of course here we are changing the edge orders only temporarily, and when the linear test does not show rigidity we restore the original orders and go to the next step.
\end{remark}

\begin{example} \label{exp:cube}
Figure 3 shows a compact 3-dimensional hyperbolic cube \textsf{cu23}.

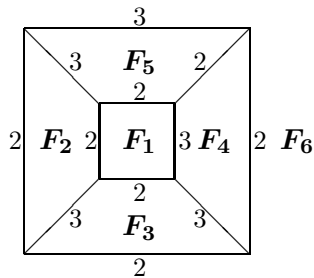
\begin{figure}[h]
\centering
\setlength{\unitlength}{1cm}
\begin{picture}(3,3.5)
\put(1.3,1.4){\large $\boldsymbol{F_1}$}
\put(0.2,1.4){\large $\boldsymbol{F_2}$}
\put(1.3,0.25){\large $\boldsymbol{F_3}$}
\put(2.3,1.4){\large $\boldsymbol{F_4}$}
\put(1.3,2.45){\large $\boldsymbol{F_5}$}
\put(3.4,1.4){\large $\boldsymbol{F_6}$}
    \put(1.45,2.05){$2$}
    \put(1,2){\line(1,0){1}}
    \put(2.05,1.4){$3$}
    \put(2,2){\line(0,-1){1}}
    \put(1.45,0.7){$2$}
    \put(2,1){\line(-1,0){1}}
    \put(0.8,1.4){$2$}
    \put(1,1){\line(0,1){1}}
    \put(0.6,2.45){$3$}
    \put(1,2){\line(-1,1){1}}
    \put(2.25,2.45){$2$}
    \put(2,2){\line(1,1){1}}
    \put(2.25,0.35){$3$}
    \put(2,1){\line(1,-1){1}}
    \put(0.6,0.35){$3$}
    \put(1,1){\line(-1,-1){1}}
    \put(-0.2,1.4){$2$}
    \put(0,0){\line(0,1){3}}
    \put(1.45,3.05){$3$}
    \put(0,3){\line(1,0){3}}
    \put(3.05,1.4){$2$}
    \put(3,3){\line(0,-1){3}}
    \put(1.45,-0.3){$2$}
    \put(3,0){\line(-1,0){3}}
\end{picture}
\caption{Labels of faces and edges of \textsf{cu23}}
\end{figure}
\noindent
Then the sets of rigid faces at level 1, level 2, and level 3 are $\{ F_1, F_6 \}$, $\{ F_4, F_5 \}$, and $\{F_2, F_3\}$, respectively. Hence every face of \textsf{cu23} occurs as a rigid face at level $1$, $2$ or $3$, and thus \textsf{cu23} is projectively rigid relative to the mirrors.
\end{example}

\subsection{Notations for figures and tables}\label{subsec:notation}

The following notations will be used in the figures and tables throughout this paper:
\begin{itemize}
\item Each $e_i$ is an edge order, corresponding to a dihedral angle $\pi/e_i$,
\item O = the number of variables $-$ the number of Vinberg equations (\ref{eqn:normalize})--(\ref{eqn:Vinberg2}),
\item I = the dimension of infinitesimal restricted deformation space of real projective structures,
\item A = the dimension of local restricted deformation space of real projective structures,
\item L = Is it possible to apply the linear test of rigidity? (yes or no), and the maximum level needed,
\item  J = Does the calculation of the Jacobian $D$ give a full description of the local restricted deformation space? (yes or no),
\item S = the minimum of the singular values of the Jacobian $D$.
\end{itemize}

\subsection{The results for cubes}\label{subsec:cubes}

Let $P$ be a compact hyperbolic cube, all of whose dihedral angles are $\pi/2$ or $\pi/3$.
By step 1 in \textsection \ref{subsec:main_alg}, the total number of such cubes is $34$ (up to symmetries).
See Table 1.
These orbifolds were tabulated by using Matlab to check the conditions of Andreev's theorem; the Matlab files used are available from the web page  (see cu.m in \cite{Lee}.)

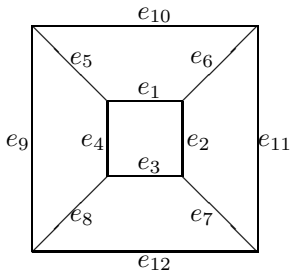
\begin{figure}[h]
\centering
\setlength{\unitlength}{1cm}
\begin{picture}(3,3.5)

    \put(1.4,2.1){$e_1$}
    \put(1,2){\line(1,0){1}}
    \put(2.05,1.4){$e_2$}
    \put(2,2){\line(0,-1){1}}
    \put(1.4,1.1){$e_3$}
    \put(2,1){\line(-1,0){1}}
    \put(0.65,1.4){$e_4$}
    \put(1,1){\line(0,1){1}}

    \put(0.5,2.5){$e_5$}
    \put(1,2){\line(-1,1){1}}
    \put(2.1,2.5){$e_6$}
    \put(2,2){\line(1,1){1}}
    \put(2.1,0.4){$e_7$}
    \put(2,1){\line(1,-1){1}}
    \put(0.5,0.4){$e_8$}
    \put(1,1){\line(-1,-1){1}}

    \put(-0.35,1.4){$e_9$}
    \put(0,0){\line(0,1){3}}
    \put(1.4,3.1){$e_{10}$}
    \put(0,3){\line(1,0){3}}
    \put(3.0,1.4){$e_{11}$}
    \put(3,3){\line(0,-1){3}}
    \put(1.4,-0.2){$e_{12}$}
    \put(3,0){\line(-1,0){3}}
\end{picture}
\caption{Labels of edges of cubes}
\end{figure}
\begin{table}
\centering
\begin{tabular}{|c|c|c|c|c||c|c|}
\hline
name & $e_1 e_2 \cdots e_{11} e_{12}$ & O & I & A & L &  J \\
\hline
\textsf{cu1} & $2 3 2 2 2 2 2 3 2 2 2 3$ & -3 & 0 & 0 & yes, level 2 & $\cdot$ \\ 
\hline
\textsf{cu2} & $2 3 2 2 2 2 2 3 2 2 3 3$ & -2 & 0 & 0 & yes, level 3 & $\cdot$ \\ 
\hline
\textsf{cu3} & $2 3 2 2 2 2 2 3 2 3 2 2$ & -3 & 0 & 0 & yes, level 1 & $\cdot$  \\ 
\hline
\textsf{cu4} & $2 3 2 2 2 2 2 3 2 3 2 3$ & -2 & 0 & 0 & yes, level 2 & $\cdot$ \\ 
\hline
\textsf{cu5} & $2 3 2 2 2 2 2 3 2 3 3 3$ & -1 & 0 & 0 & yes, level 3 & $\cdot$ \\ 
\hline
\textsf{cu6} & $2 3 2 2 2 2 2 3 3 3 2 2$ & -2 & 0 & 0 & yes, level 2 & $\cdot$ \\ 
\hline
\textsf{cu7} & $2 3 2 2 2 2 2 3 3 3 3 2$ & -1 & 0 & 0 & yes, level 3 & $\cdot$ \\ 
\hline
\textsf{cu8} & $2 3 2 2 2 2 3 2 2 2 2 3$ & -3 & 0 & 0 & yes, level 2 & $\cdot$ \\ 
\hline
\textsf{cu9} & $2 3 2 2 2 2 3 2 2 3 3 2$ & -2 & 0 & 0 & yes, level 2 & $\cdot$ \\ 
\hline
\textsf{cu10} & $2 3 2 2 2 2 3 2 3 2 2 3$ & -2 & 0 & 0 & yes, level 3 & $\cdot$ \\ 
\hline
\textsf{cu11} & $2 3 2 2 2 2 3 2 3 3 2 2$ & -2 & 0 & 0 & yes, level 2 & $\cdot$ \\ 
\hline
\textsf{cu12} & $2 3 2 2 2 2 3 2 3 3 2 3$ & -1 & 0 & 0 & yes, level 3 & $\cdot$ \\ 
\hline
\textsf{cu13} & $2 3 2 2 2 2 3 2 3 3 3 2$ & -1 & 0 & 0 & yes, level 2 & $\cdot$ \\ 
\hline
\textsf{cu14} & $2 3 2 2 2 2 3 3 3 3 2 2$ & -1 & 0 & 0 & yes, level 3 & $\cdot$ \\ 
\hline
\textsf{cu15} & $2 3 2 2 2 2 3 3 3 3 3 2$ & 0 & 0 & 0 & no & yes
\\ 
\hline
\textsf{cu16} & $2 3 2 2 2 3 2 3 3 3 2 2$ & -1 & 0 & 0 & yes, level 3 & $\cdot$ \\ 
\hline
\textsf{cu17} & $2 3 2 2 2 3 3 2 2 3 2 3$ & -1 & 1 & \textbf{1} & no & no \\ 
\hline
\textsf{cu18} & $2 3 2 2 2 3 3 2 3 3 2 3$ & 0 & 1 & \textbf{1} & no & no
\\ 
\hline
\textsf{cu19} & $2 3 2 2 2 3 3 3 3 3 2 2$ & 0 & 0 & 0 & no & yes
\\ 
\hline
\textsf{cu20} & $2 3 2 2 3 2 2 3 2 2 3 3$ & -1 & 0 & 0 & yes, level 3 & $\cdot$  \\ 
\hline
\textsf{cu21} & $2 3 2 2 3 2 2 3 2 3 2 3$ & -1 & 1 & \textbf{1} & no & no \\ 
\hline
\textsf{cu22} & $2 3 2 2 3 2 2 3 2 3 3 3$ & 0 & 1 & \textbf{1} & no & no  \\ 
\hline
\textsf{cu23} & $2 3 2 2 3 2 3 3 2 3 2 2$ & -1 & 0 & 0 & yes, level 3 & $\cdot$ \\ 
\hline
\textsf{cu24} & $2 3 2 2 3 2 3 3 2 3 2 3$ & 0 & 0 & 0 & no & yes
\\ 
\hline
\textsf{cu25} & $2 3 2 2 3 2 3 3 2 3 3 2$ & 0 & 0 & 0 & no & yes
\\ 
\hline
\textsf{cu26} & $2 3 2 2 3 3 3 3 2 2 2 3$ & 0 & 1 & \textbf{0} & no & no
\\ 
\hline
\textsf{cu27} & $2 3 2 2 3 3 3 3 2 3 2 3$ & 1 & 2 & \textbf{1} & no & no
\\ 
\hline
\textsf{cu28} & $2 3 2 3 2 2 2 3 2 2 3 3$ & -1 & 0 & 0 & no & yes
\\ 
\hline
\textsf{cu29} & $2 3 2 3 2 3 2 3 2 3 2 3$ & 0 & 1 & \textbf{0} & no & no
\\ 
\hline
\textsf{cu30} & $2 3 2 3 2 3 3 2 3 3 2 3$ & 1 & 1 & 1 & no & yes
\\ 
\hline
\textsf{cu31} & $2 3 2 3 2 3 3 3 2 3 2 3$ & 1 & 1 & 1 & no & yes
\\ 
\hline
\textsf{cu32} & $2 3 2 3 2 3 3 3 3 3 2 2$ & 1 & 1 & 1 & no & yes
\\ 
\hline
\textsf{cu33} & $2 3 2 3 3 3 3 3 2 3 2 3$ & 2 & 3 & \textbf{2} & no & no
\\ 
\hline
\textsf{cu34} & $2 3 3 2 2 3 2 3 3 3 2 2$ & 0 & 1 & \textbf{1} & no & no
\\ 
\hline
\end{tabular}
\caption{The list of cubes. (The bold letters mean that the results are obtained from the final step of the main algorithm of \textsection 3.1.)}
\end{table}

By step 2 in \textsection \ref{subsec:main_alg} (the linear test for rigidity), we find that seventeen cubes (\textsf{cu1-cu14, cu16, cu20, cu23}) are projectively rigid relative to the mirrors.

Using steps 3-5 in \textsection \ref{subsec:main_alg}, exact algebraic computations of the dimensions of Zariski tangent spaces determine the dimensions of the local restricted deformation spaces for eight cubes (\textsf{cu15, cu19, cu24, cu25, cu28, cu30-cu32})
since we can show that the Jacobian matrices $D$ for these cubes are of full rank. The computations were done with arbitrary precision
algorithms in Mathematica, and were also numerically verified using Matlab.

We need to look beyond the Jacobian matrix to calculate the dimension of the local restricted deformation spaces for the remaining cubes
(\textsf{cu17, cu18, cu21, cu22, cu26, cu27, cu29, cu33, cu34}) as the Jacobian matrices of these cubes are rank-deficient.
Using step 6 in \textsection \ref{subsec:main_alg}, we instead obtain a Gr\"{o}bner basis for the ideal generated by Vinberg's equations (\ref{eqn:normalize})--(\ref{eqn:Vinberg2})
using Mathematica.
The detailed computations for \textsf{cu\,$i$} is saved in a Mathematica notebook file (\textsf{cu\,$i$}.nb) available from the web page \cite{Lee} for each $i=15,17,18,19,21,22,24$,
$25,26,27,28,29,30,31,32,33,34$.
(Note that choosing the hyperbolic
solutions for the cubes  is often a nontrivial process involving geometric considerations, and the choices are explained in the files themselves.)

\begin{example} \label{example:cu21}
As an example, we use \textsf{cu21} to illustrate the method for computing the local restricted deformation space. (See Figure 5.)
All the local restricted deformation spaces for the cubes except \textsf{cu27} and \textsf{cu33} are obtained by a similar method.

\begin{figure}[h]
\centering
\setlength{\unitlength}{1cm}
\begin{picture}(3,3.5)
\put(1.3,1.35){\large $\boldsymbol{F_2}$}
\put(0.2,1.35){\large $\boldsymbol{F_5}$}
\put(1.3,0.25){\large $\boldsymbol{F_6}$}
\put(2.3,1.35){\large $\boldsymbol{F_3}$}
\put(1.3,2.45){\large $\boldsymbol{F_1}$}
\put(3.4,1.35){\large $\boldsymbol{F_4}$}
    \put(1.45,2.05){$2$}
    \put(1,2){\line(1,0){1}}
    \put(2.05,1.4){$3$}
    \put(2,2){\line(0,-1){1}}
    \put(1.45,0.7){$2$}
    \put(2,1){\line(-1,0){1}}
    \put(0.8,1.4){$2$}
    \put(1,1){\line(0,1){1}}
    \put(0.6,2.45){$3$}
    \put(1,2){\line(-1,1){1}}
    \put(2.25,2.45){$2$}
    \put(2,2){\line(1,1){1}}
    \put(2.25,0.35){$2$}
    \put(2,1){\line(1,-1){1}}
    \put(0.6,0.35){$3$}
    \put(1,1){\line(-1,-1){1}}
    \put(-0.2,1.4){$2$}
    \put(0,0){\line(0,1){3}}
    \put(1.45,3.05){$3$}
    \put(0,3){\line(1,0){3}}
    \put(3.05,1.4){$2$}
    \put(3,3){\line(0,-1){3}}
    \put(1.45,-0.3){$3$}
    \put(3,0){\line(-1,0){3}}
\end{picture}
\caption{\textsf{cu21}}
\end{figure}

First, we note that the set of all rigid faces is $\{F_2, F_3\}$, and these faces are at level $1$. So the linear test of rigidity is not applicable.

Second, we find the unit normals $\nu_i$ for \textsf{cu21} as follows.
We must solve the system $\{\Psi_k=0\}_{k=1}^n$ of hyperbolic equations (\ref{eqn:hyp}), where $n=f+e=18$.
We choose the vertex $F_{123}=F_1 \cap F_2 \cap F_3$ whose adjacent three edges have orders $(2,2,3)$.
Then, by an isometry of $\H^3$,  we can assume that the normals of the adjacent faces are
\begin{displaymath}
\nu_1=(0,1,0,0),\;\; \nu_2=(0,0,1,0),\;\; \text{and}\;\; \nu_3=(0,0,-\frac{1}{2},\frac{\sqrt{3}}{2}).
\end{displaymath}
Hence we satisfy the following six hyperbolic equations:
{\setlength\arraycolsep{1pt}
\begin{eqnarray*}
& & \langle \nu_1, \nu_2 \rangle = \langle \nu_1, \nu_3 \rangle = 0 \;\; \text{and} \;\; \langle \nu_2, \nu_3 \rangle=-1/2,\\
& & \langle \nu_1, \nu_1 \rangle = \langle \nu_2, \nu_2 \rangle = \langle \nu_3, \nu_3 \rangle = 1.
\end{eqnarray*}}

Since the orders of the two edges $F_{14}$ and $F_{34}$ are $3$ and $2$ respectively,  we let
\begin{displaymath}
\nu_4=(x,-\frac{1}{2},u,\frac{u}{\sqrt{3}}).
\end{displaymath}
Similarly, we let
\begin{displaymath}
\nu_5=(y,-\frac{1}{2},0,\frac{2v}{\sqrt{3}})\;\; \text{and} \;\; \nu_6=(z,w,0,0).
\end{displaymath}
Then we satisfy the six hyperbolic equations
{\setlength\arraycolsep{1pt}
\begin{eqnarray*}
& & \langle \nu_1, \nu_4 \rangle = \langle \nu_1, \nu_5 \rangle = -1/2,\\
& & \langle \nu_3, \nu_4 \rangle = \langle \nu_2, \nu_5 \rangle = \langle \nu_2, \nu_6 \rangle = \langle \nu_3, \nu_6 \rangle = 0.
\end{eqnarray*}}
Therefore, we must solve the remaining six hyperbolic equations
{\setlength\arraycolsep{1pt}
\begin{eqnarray*}
& & \langle \nu_4, \nu_5 \rangle =0 \;\; \text{and}\;\; \langle \nu_4, \nu_6 \rangle = \langle \nu_5, \nu_6 \rangle = -1/2, \\
& & \langle \nu_4, \nu_4 \rangle = \langle \nu_5, \nu_5 \rangle = \langle \nu_6, \nu_6 \rangle = 1.
\end{eqnarray*}}
However, these equations have many solutions.
Among them we choose a solution such that non-diagonal entries of the Gram matrix $G=(g_{ij})$ with $g_{ij}=\langle \nu_i, \nu_j \rangle$ are non-positive.
In particular, $u$, $v$ and $w$ are non-positive since $g_{24}=u$, $g_{35}=v$ and $g_{16}= w$.
This leaves two solutions:
{\setlength\arraycolsep{1pt}
\begin{eqnarray}
x=y=z=-\frac{\sqrt{5}}{2},\;\; u=v=-\sqrt{\frac{3}{2}},\;\; \text{and}\;\; w=-\frac{3}{2} \label{soln1}\\
x=y=z=\frac{\sqrt{5}}{2},\;\; u=v=-\sqrt{\frac{3}{2}},\;\; \text{and}\;\; w=-\frac{3}{2} \label{soln2}
\end{eqnarray}}
For solution (\ref{soln1}), the cone defined by the inequalities $\langle X, \nu_i \rangle \geq 0$ intersects $\mathbb{H}^3$ since $X=(1,0,0,0)$ satisfies the inequalities $\langle X, \nu_i \rangle = -\nu_{i1} \geq 0$.
However, in the case of solution (\ref{soln2}), the cone defined by the inequalities $\langle X, \nu_i \rangle \leq 0$ intersects $\mathbb{H}^3$ since $X=(1,0,0,0)$ satisfies the inequalities $\langle X, \nu_i \rangle = -\nu_{i1} \leq 0$.
Here, the cone $K$ is defined by the inequalities $\langle X, \nu_i \rangle \leq 0$.
Hence only solution (\ref{soln2}) is appropriate for our hyperbolic cube.
Furthermore, the Gram matrix $G$ is as follows:
\begin{displaymath}
G =
\left( \begin{array}{cccccc}
1 & 0 & 0 & -\frac{1}{2} & -\frac{1}{2} & -\frac{3}{2} \\
0 & 1 & -\frac{1}{2} & -\sqrt{\frac{3}{2}} & 0 & 0\\
0 & -\frac{1}{2} & 1 & 0 & -\sqrt{\frac{3}{2}} & 0 \\
-\frac{1}{2} & -\sqrt{\frac{3}{2}} & 0 & 1 & 0 & -\frac{1}{2} \\
-\frac{1}{2} & 0 & -\sqrt{\frac{3}{2}} & 0 & 1 & -\frac{1}{2} \\
-\frac{3}{2} & 0 & 0 & -\frac{1}{2} & -\frac{1}{2} & 1
\end{array} \right).
\end{displaymath}
Hence the solution (\ref{soln2}) is the unique solution satisfying
the conditions of Theorem \ref{thm:cube}, and
the cone $K$ gives the 3-dimensional compact hyperbolic cube \textsf{cu21}$=K\cap \mathbb{H}^3$.

\begin{remark}\label{rmk:cube2}
For each cube,  we start by choosing a  vertex whose adjacent edges have orders $(2,2,3)$ or $(2,2,2)$.
If the edge orders are $2,2,3$ we choose normals $\nu_1, \nu_2, \nu_3$ for the adjacent three faces as above.
This applies to the cases  \textsf{cu21, cu22, cu24, cu25, cu26, cu29, cu30, cu31} and \textsf{cu32}.
If the edge orders are $2,2,2$  we let
 \begin{displaymath}
\nu_1=(0,1,0,0),\;\; \nu_2=(0,0,1,0)\;\; \text{and},\;\; \nu_3=(0,0,0,1).
\end{displaymath}
This applies to the cases \textsf{cu15, cu17, cu18, cu19, cu28} and \textsf{cu34}.
\end{remark}

Third, using the linear functionals $\alpha_i=J\nu_i$ we form the Jacobian matrix $D=D\Phi$ for Vinberg's equations (\ref{eqn:normalize})--(\ref{eqn:Vinberg2}) at the hyperbolic point.
Note that $D$ is a $25 \times 24$ matrix.
Using Mathematica we find that the rank of $D$ is $23$, and so $D$ is rank-deficient.
Since the dimension of kernel of the Jacobian matrix $D$ is $1$,
the dimension of the infinitesimal restricted deformation space of real projective structures is $1$.

Finally, to obtain the dimension of the local restricted deformation space of \textsf{cu21},
we compute a Gr\"{o}bner basis of the ideal $\langle \Phi_1, \ldots, \Phi_N \rangle$ with $N=f+e+e_2=25$.
Before doing this, we introduce new coordinates on $\mathbb{R}^{4f}=\mathbb{R}^{24}$ by letting $c_i=b_i-t_i$, where $t=\{t_i\}$
corresponds to the hyperbolic point in $\mathfrak{D}_P(\hat{P})$.
(Relative to this coordinate system, $t$ is the origin.)
We compute a Gr\"{o}bner basis of the ideal $\langle \Phi_1, \ldots, \Phi_{25} \rangle$ with respect to the lexicographic order with
$
c_{41}<c_{42}<c_{43}<c_{44} < c_{51}< c_{52} < c_{53} < c_{54} < c_{61} < c_{62} < c_{63} < c_{64} < c_{11} < c_{12} < c_{13} < c_{14}
 < c_{2 1} < c_{22} < c_{23} < c_{24} < c_{31} < c_{32} < c_{33} < c_{34}.
$
Then the Gr\"{o}bner basis for \textsf{cu21} is
{\setlength\arraycolsep{1pt}
\begin{eqnarray*}
& & \{ c_{34}, c_{33}, c_{32}, c_{31}, c_{24}, c_{23}, c_{22}, c_{21}, c_{14}, c_{13}, c_{12}, c_{64}, c_{63}, \\
& & \mathbf{-c_{11}+\frac{2c_{62}}{\sqrt{5}}+2c_{11}c_{62}}, \sqrt{5}c_{61}+3c_{62}, c_{54}, c_{53}, \mathbf{-c_{11}+\frac{2c_{52}}{\sqrt{5}}+c_{11}c_{52}}, \\
& & \mathbf{-c_{52}+c_{62}+c_{52}c_{62}}, \sqrt{5}c_{51}+c_{52}, c_{44}, c_{43}, c_{42}-c_{52}, \sqrt{5}c_{41}+c_{52}\}
\end{eqnarray*}}
Let $f_i$ be the $i$th
polynomial shown in bold letters.
The Gr\"{o}bner basis implies that
\begin{itemize}
\item $c_{11}$ is a free variable, and $c_{11}$ determines $c_{62}$ and $c_{52}$,
\item $c_{62}$ determines $c_{61}$,  and $c_{52}$ determines $c_{51}$, $c_{42}$ and $c_{41}$,
\item the remaining variables are zero, and
\item the equation $f_3=0$ is implied by the equations $f_1=0$ and $f_2=0$.

\end{itemize}
The last observation follows since the three polynomials in bold letters satisfy the
relation
\begin{displaymath}
c_{62}f_2+c_{11}f_3=c_{52}f_1.
\end{displaymath}
In other words, if $c_{11}\neq 0$ then $f_3=(c_{52}f_1-c_{62}f_2)/c_{11}$. Hence $f_1=f_2=0$ implies $f_3=0$.
Also, if $c_{11}=0$ then $f_1=f_2=0$ imply $c_{62}=c_{52}=0$, and thus $f_3=0$.
Therefore, this implies that the dimension of the local restricted deformation space is also 1.
\end{example}

\begin{example} \label{example:cu27}
As another example, we describe the calculation of the local restricted deformation space for \textsf{cu27},
noting the differences to the method used for \textsf{cu21}. (See Figure 6).

\begin{figure}[h]
\centering
\setlength{\unitlength}{1cm}
\begin{picture}(3,3.5)
\put(1.3,1.35){\large $\boldsymbol{F_1}$}
\put(0.2,1.35){\large $\boldsymbol{F_3}$}
\put(1.3,0.25){\large $\boldsymbol{F_4}$}
\put(2.3,1.35){\large $\boldsymbol{F_5}$}
\put(1.3,2.45){\large $\boldsymbol{F_2}$}
\put(3.4,1.35){\large $\boldsymbol{F_6}$}
    \put(1.45,2.05){$2$}
    \put(1,2){\line(1,0){1}}
    \put(2.05,1.4){$3$}
    \put(2,2){\line(0,-1){1}}
    \put(1.45,0.7){$2$}
    \put(2,1){\line(-1,0){1}}
    \put(0.8,1.4){$2$}
    \put(1,1){\line(0,1){1}}
    \put(0.6,2.45){$3$}
    \put(1,2){\line(-1,1){1}}
    \put(2.25,2.45){$3$}
    \put(2,2){\line(1,1){1}}
    \put(2.25,0.35){$3$}
    \put(2,1){\line(1,-1){1}}
    \put(0.6,0.35){$3$}
    \put(1,1){\line(-1,-1){1}}
    \put(-0.2,1.4){$2$}
    \put(0,0){\line(0,1){3}}
    \put(1.45,3.05){$3$}
    \put(0,3){\line(1,0){3}}
    \put(3.05,1.4){$2$}
    \put(3,3){\line(0,-1){3}}
    \put(1.45,-0.3){$3$}
    \put(3,0){\line(-1,0){3}}
\end{picture}
\caption{\textsf{cu27}}
\end{figure}
Steps 1-2 are similar to those for \textsf{cu21}. (The details are omitted.)
To find the exact values of $\nu_i$ using Mathematica,
we use the reflectional symmetry of \textsf{cu27}, interchanging $F_2$ and $F_4$.
Hence we let
{\setlength\arraycolsep{1pt}
\begin{eqnarray*}
& & \nu_1 = (0,0,0,1),\;\; \nu_2 = (x,-1/2,-u,0),\;\; \nu_3 = (0,1,0,0), \\
& & \nu_4 = (x,-1/2,u,0), \;\; \nu_5 = (y,v,0,-1/2), \;\;\nu_6 = (z,0,0,w).
\end{eqnarray*}}
By similar computations to those for \textsf{cu21}, we obtain the following:
\begin{displaymath}
\begin{array}{|c|c|c|}
\hline
\text{notation} & \text{numerical value} & \text{a real root of} \\
\hline
u & -1.36278 & 5 + 44 t^2 - 144 t^4 + 64 t^6=0 \\
\hline
v & -1.45161 & -1 - 2 t + 2 t^2 + 2 t^3=0 \\
\hline
w & -1.10716 & 1 - 4 t + 4 t^3=0 \\
\hline
x & 1.05222 & -1 - 4 t^2 + 4 t^6=0 \\
\hline
y & 1.16497 & -1 - 52 t^2 - 48 t^4 + 64 t^6=0 \\
\hline
z & 0.47519 & -1 + 16 t^4 + 16 t^6=0 \\
\hline
\end{array}
\end{displaymath}

We again compute the rank of the Jacobian matrix $D$ at the hyperbolic point.
Here $D$ is a $23 \times 24$ matrix with rank $22$. Thus $D$ does not have full rank,
and the dimension of kernel of the Jacobian matrix $D$ is $2$.

Since the expressions of $\alpha_i=J\nu_i$ for \textsf{cu27} are complicated,
Mathematica is unable to compute a Gr\"{o}bner basis of the ideal $\langle \Phi_1, \ldots, \Phi_{23} \rangle$ in a reasonable time.
To make the problem easier for Mathematica to solve,
we find an algebraic number $\theta$ such that $\{u,v,w,x,y,z\} \subset \mathbb{Q}(\theta)$. Here
$ \theta \approx -0.395609$ and is a real root of $64 - 384 t^2 - 208 t^4 + 320 t^6 - 52 t^8 - 24 t^{10} + t^{12}=0$,
and we convert the above $u$, $v$, $w$, $x$, $y$, and $z$ to elements of $\mathbb{Q}(\theta)$ as follows:
\begin{displaymath}
\begin{array}{|c|c|}
\hline
\text{notation} & \text{as a element of } \mathbb{Q}(\theta) \\
\hline
u & \frac{13}{4}\theta+\frac{13}{8}\theta^3-\frac{5}{2}\theta^5+\frac{13}{32}\theta^7+\frac{3}{16}\theta^9-\frac{1}{128}\theta^{11}  \\
\hline
v & -\frac{47}{32}-\frac{1}{16}\theta^2+\frac{73}{64}\theta^4-\frac{1}{4}\theta^6-\frac{47}{512}\theta^8+\frac{1}{256}\theta^{10}  \\
\hline
w & -1-\frac{7}{8}\theta^2+\frac{5}{4}\theta^4-\frac{13}{64}\theta^6-\frac{3}{32}\theta^8+\frac{1}{256}\theta^{10}  \\
\hline
x & -\frac{169}{64}\theta-\frac{45}{128}\theta^3+\frac{195}{128}\theta^5-\frac{89}{256}\theta^7-\frac{129}{1024}\theta^9+\frac{11}{2048}\theta^{11}  \\
\hline
y & -\frac{11}{4}\theta-\frac{13}{8}\theta^3+\frac{5}{2}\theta^5-\frac{13}{32}\theta^7-\frac{3}{16}\theta^9+\frac{1}{128}\theta^{11} \\
\hline
z & -\frac{35}{32}\theta-\frac{55}{64}\theta^3+\frac{9}{8}\theta^5-\frac{7}{64}\theta^7-\frac{37}{512}\theta^9+\frac{3}{1024}\theta^{11} \\
\hline
\end{array}
\end{displaymath}

Next, we introduce new coordinates on $\mathbb{R}^{4f}=\mathbb{R}^{24}$ by letting $c_i=b_i-t_i$, where $t=\{t_i\}$
correspond to the hyperbolic point in $\mathfrak{D}_P(\hat{P})$.
We compute a Gr\"{o}bner basis of the ideal $\langle \Phi_1, \ldots, \Phi_{23} \rangle$ with respect to the lexicographic order with
$
c_{21} <  c_{22} < c_{23} < c_{24}  < c_{41} < c_{42} < c_{43} < c_{44} < c_{51} < c_{52} < c_{53} < c_{54} < c_{61} < c_{62} < c_{63} <
c_{64} < c_{31} < c_{32} < c_{33} < c_{34} < c_{11} < c_{12} < c_{13} < c_{14}.
$
Then the Gr\"{o}bner basis for \textsf{cu27} is
\begin{displaymath}
\{ c_{14}, c_{13}, c_{12}, c_{11}, c_{34}, c_{32}, c_{31}, c_{64}, (c_{33}-c_{63})^2, c_{62}, c_{61}, c_{54}, \ldots \}.
\end{displaymath}
Replacing the polynomial $(c_{33}-c_{63})^2$ by $(c_{33}-c_{63})$ gives a new ideal
$\langle  c_{33}-c_{63}, \Phi_1, \ldots, \Phi_{23} \rangle$ with the same underlying solution set.
This has  the Gr\"{o}bner basis
{\setlength\arraycolsep{1pt}
\begin{eqnarray*}
\{ & & c_{14}, c_{13}, c_{12}, c_{11}, c_{34}, c_{32}, c_{31}, c_{64}, c_{33}-c_{63}, c_{62}, c_{61}, c_{54}, c_{33}-c_{53},\\
& &  c_{52}, c_{51}, c_{44}, \mathbf{-\delta c_{33}+\varepsilon c_{43}+c_{33}c_{43}}, -\delta c_{42}+c_{43}, \mu c_{41}+c_{43}, c_{24},\\
& &  \mathbf{-\varepsilon c_{23}+\delta c_{33}+c_{23}c_{33}}, \mathbf{-\delta c_{23}+\delta c_{43}+2c_{23}c_{43}}, \delta c_{22}+c_{23}, -\mu c_{21}+c_{23}\},
\end{eqnarray*}}
where
\begin{displaymath}
\begin{array}{|c|c|c|}
\hline
\text{notation} & \text{numerical value} & \text{a real root of} \\
\hline
\delta & 1.25796 & 961 - 1017 t^2 + 251 t^4 + 5 t^6=0 \\
\hline
\varepsilon & 0.733796 &  64 - 144 t^2 + 44 t^4 + 5 t^6=0 \\
\hline
\mu & 0.597767 & 961 - 2288 t^2 - 1152 t^4 + 80 t^6=0 \\
\hline
\end{array}
\end{displaymath}
Let $f_i$ be the $i$th
polynomial shown in bold letters.
Then this Gr\"{o}bner basis implies that
\begin{itemize}
\item $c_{33}$ is a free variable, and $c_{33}$ determines $c_{63}$, $c_{53}$, $c_{43}$ and $c_{23}$,
\item $c_{43}$ determines $c_{42}$ and $c_{41}$, and $c_{23}$ determines $c_{22}$ and $c_{21}$,
\item the remaining variables are zeros, and
\item the equation $f_3=0$ is implied by the equations $f_1=0$ and $f_2=0$
where $c_{2,3} f_1 + c_{4,3} f_2 = c_{3,3} f_3$ holds.

\end{itemize}
Thus, although the dimension of the infinitesimal restricted deformation space is 2, the dimension of the local restricted deformation space is 1.
The details of the calculations are available from the webpage \cite{Lee}.

\end{example}

\subsection{The results for dodecahedra}\label{subsec:dodec}

Let $P$ be a $3$-dimensional compact hyperbolic dodecahedron, all of whose dihedral angles are $\pi/2$ or $\pi/3$.
We assume that each face has less than three edges of order $2$.
Then the total number of such orbifolds is $13$ up to symmetries (see table 2).
These results were obtained by using
Matlab to check the conditions of Andreev's theorem (see do.m in \cite{Lee}).

\begin{figure}[h]
\centering
\setlength{\unitlength}{1cm}
\begin{picture}(5,5.5)
\put(2.7,3.0){$e_1$}
\put(2.85,2.45){$e_2$}
\put(2.35,2.05){$e_3$}
\put(1.85,2.45){$e_4$}
\put(2.05,3.0){$e_5$}
\put(2.5,3.65){$e_6$}
\put(3.35,3.05){$e_7$}
\put(3.15,1.8){$e_8$}
\put(1.5,1.8){$e_9$}
\put(1.3,3.05){$e_{10}$}
\put(2.75,4.15){$e_{11}$}
\put(4.0,2.65){$e_{12}$}
\put(2.9,1.1){$e_{13}$}
\put(0.8,1.7){$e_{14}$}
\put(0.85,3.6){$e_{15}$}
\put(3.7,3.6){$e_{16}$}
\put(3.7,1.65){$e_{17}$}
\put(1.7,1.05){$e_{18}$}
\put(0.55,2.65){$e_{19}$}
\put(1.8,4.15){$e_{20}$}
\put(3.3,4.4){$e_{21}$}
\put(4.3,2.15){$e_{22}$}
\put(2.5,0.6){$e_{23}$}
\put(0.25,2.2){$e_{24}$}
\put(1.2,4.4){$e_{25}$}
\put(2.3,4.8){$e_{26}$}
\put(4.4,3.4){$e_{27}$}
\put(3.6,0.9){$e_{28}$}
\put(0.95,0.9){$e_{29}$}
\put(0.1,3.4){$e_{30}$}
\put(2,2){\line(1,0){1}}
\put(2,2){\line(-1,3){0.3}}
\put(3,2){\line(1,3){0.3}}
\put(2.5,3.5){\line(4,-3){0.8}}
\put(2.5,3.5){\line(-4,-3){0.8}}
\put(2.5,3.5){\line(0,1){0.5}}
\put(2.5,4){\line(6,1){1.1}}
\put(2.5,4){\line(-6,1){1.1}}
\put(3.3,2.9){\line(3,1){0.5}}
\put(1.7,2.9){\line(-3,1){0.5}}
\put(3.8,3.07){\line(-1,6){0.185}}
\put(1.2,3.07){\line(1,6){0.185}}
\put(3,2){\line(3,-4){0.34}}
\put(2,2){\line(-3,-4){0.34}}
\put(3.8,3.07){\line(1,-2){0.46}}
\put(1.2,3.07){\line(-1,-2){0.46}}
\put(3.35,1.55){\line(3,2){0.90}}
\put(1.65,1.55){\line(-3,2){0.90}}
\put(2.5,0.87){\line(-5,4){0.84}}
\put(2.5,0.87){\line(5,4){0.84}}
\put(2.5,0.87){\line(0,-1){0.69}}
\put(1,4.7){\line(1,0){3}}
\put(1,4.7){\line(3,-4){0.38}}
\put(4,4.7){\line(-3,-4){0.38}}
\put(1,4.7){\line(-1,-3){0.9}}
\put(4,4.7){\line(1,-3){0.9}}
\put(0.75,2.15){\line(-4,-1){0.65}}
\put(4.25,2.15){\line(4,-1){0.65}}
\put(2.5,0.18){\line(-4,3){2.4}}
\put(2.5,0.18){\line(4,3){2.4}}
\end{picture}
\caption{Labels of edges of a dodecahedron}
\end{figure}
\begin{table}[h]
\centering
\begin{tabular}{|c|c|c|c|c||c|c|}
\hline
name & $e_1\:e_2 \cdots e_{29}\:e_{30}$ & O & I & A & J & S \\
\hline
\textsf{do1} & $2 3 2 3 3 2 3 2 3 2 3 2 3 2 3 3 3 3 3 3 2 3 2 3 2 3 3 2 3 2$ & -6 & 0 & 0 & yes & 0.17653\\
\hline
\textsf{do2} & $2 3 2 3 3 2 3 3 3 2 3 2 2 2 3 3 3 3 3 3 2 3 3 2 2 3 3 2 3 3$ & -5 & 0 & 0 & yes & 0.13121\\
\hline
\textsf{do3} & $2 3 2 3 3 2 3 3 3 2 3 2 2 2 3 3 3 3 3 3 2 3 3 3 2 3 3 2 2 3$ & -5 & 0 & 0 & yes & 0.14468\\
\hline
\textsf{do4} & $2 3 2 3 3 2 3 3 3 2 3 2 2 2 3 3 3 3 3 3 2 3 3 3 2 3 3 2 3 2$ & -5 & 0 & 0 & yes & 0.13707\\
\hline
\textsf{do5} & $2 3 2 3 3 2 3 3 3 2 3 2 3 2 3 3 2 2 3 3 2 3 3 3 2 3 3 2 3 2$ & -6 & 0 & 0 & yes & 0.18151\\
\hline
\textsf{do6} & $2 3 2 3 3 2 3 3 3 2 3 2 3 3 3 3 2 2 3 3 2 3 3 2 2 3 3 2 3 3$ & -5 & 0 & 0 & yes & 0.11944\\
\hline
\textsf{do7} & $2 3 2 3 3 3 2 2 3 2 3 3 3 2 3 3 3 3 3 2 2 2 2 3 3 3 3 3 3 2$ & -5 & 0 & 0 & yes & 0.12703\\
\hline
\textsf{do8} & $2 3 2 3 3 3 2 3 3 2 3 2 2 2 3 3 3 3 3 2 2 3 3 2 3 3 3 2 3 2$ & -6 & 0 & 0 & yes & 0.09580\\
\hline
\textsf{do9} & $2 3 2 3 3 3 2 3 3 2 3 2 2 2 3 3 3 3 3 2 2 3 3 3 3 3 3 2 3 2$ & -5 & 0 & 0 & yes & 0.09365\\
\hline
\textsf{do10} & $2 3 2 3 3 3 2 3 3 2 3 3 2 2 3 3 3 3 3 2 2 2 3 3 3 3 3 3 2 2$ & -5 & 0 & 0 & yes & 0.08277\\
\hline
\textsf{do11} & $2 3 2 3 3 3 3 3 3 2 2 2 2 2 3 3 3 3 3 3 3 3 3 3 2 3 2 3 2 3$ & -4 & 0 & 0 & yes & 0.06115\\
\hline
\textsf{do12} & $2 3 2 3 3 3 3 3 3 2 2 2 3 3 2 3 2 2 2 3 3 2 3 3 3 2 3 3 2 3$ & -6 & 0 & 0 & yes & 0.12412\\
\hline
\textsf{do13} & $2 3 2 3 3 3 3 3 3 2 3 3 2 2 3 2 3 3 3 2 3 2 3 3 3 2 3 3 2 3$ & -4 & 1 & \textbf{1} & no & $\cdot$ \\
\hline
\end{tabular}
\caption{The list of Dodecahedra}
\end{table}
Since we assume that each face has less than three edges of order $2$, there is no rigid face.  Hence we can skip the linear test for rigidity.

As the next step, we want to find unit normals $\nu_i$ of the dodecahedra. However, since the number of variables in the hyperbolic equations (\ref{eqn:hyp}) is large, the exact algebraic solution is hard to obtain for all cases.

\begin{example} \label{example:do13}
We first concentrate on $\textsf{do13}$, shown in Figure 8. This orbifold has rotational symmetry that will allow us to calculate its normals $\nu_i$ exactly, and also find the dimension of its local restricted deformation space.
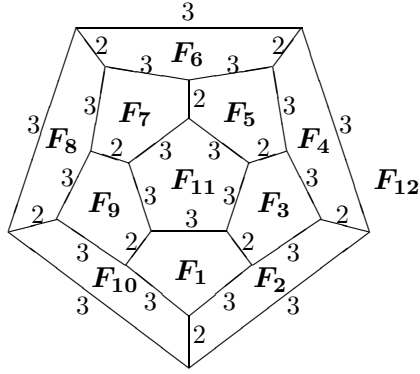
\begin{figure}
\centering
\setlength{\unitlength}{1cm}
\begin{picture}(5,5.5)
\put(2.3,1.35){\large $\boldsymbol{F_1}$}
\put(3.35,1.25){\large $\boldsymbol{F_2}$}
\put(3.40,2.25){\large $\boldsymbol{F_3}$}
\put(3.93,3.1){\large $\boldsymbol{F_{4}}$}
\put(2.95,3.45){\large $\boldsymbol{F_5}$}
\put(2.25,4.25){\large $\boldsymbol{F_{6}}$}
\put(1.55,3.45){\large $\boldsymbol{F_7}$}
\put(0.58,3.1){\large $\boldsymbol{F_{8}}$}
\put(1.15,2.25){\large $\boldsymbol{F_9}$}
\put(1.2,1.25){\large $\boldsymbol{F_{10}}$}
\put(2.25,2.55){\large $\boldsymbol{F_{11}}$}
\put(4.95,2.55){\large $\boldsymbol{F_{12}}$}
\put(2.75,2.95){$3$} 
\put(2.95,2.35){$3$}
\put(2.45,2.05){$3$}
\put(1.9,2.35){$3$}
\put(2.1,2.95){$3$}
\put(2.55,3.65){$2$} 
\put(3.4,3.0){$2$}
\put(3.2,1.75){$2$}
\put(1.65,1.75){$2$}
\put(1.45,3.0){$2$}
\put(3.0,4.10){$3$} 
\put(4.05,2.6){$3$}
\put(2.95,0.95){$3$}
\put(1.0,1.6){$3$}
\put(1.1,3.55){$3$}
\put(3.75,3.55){$3$} 
\put(3.8,1.6){$3$}
\put(1.9,0.95){$3$}
\put(0.8,2.6){$3$}
\put(1.85,4.10){$3$}
\put(3.6,4.35){$2$} 
\put(4.45,2.1){$2$}
\put(2.55,0.5){$2$}
\put(0.4,2.1){$2$}
\put(1.25,4.35){$2$}
\put(2.4,4.8){$3$} 
\put(4.5,3.3){$3$}
\put(3.8,0.9){$3$}
\put(1.0,0.9){$3$}
\put(0.35,3.3){$3$}
\put(2,2){\line(1,0){1}}
\put(2,2){\line(-1,3){0.3}}
\put(3,2){\line(1,3){0.3}}
\put(2.5,3.5){\line(4,-3){0.8}}
\put(2.5,3.5){\line(-4,-3){0.8}}
\put(2.5,3.5){\line(0,1){0.5}}
\put(2.5,4){\line(6,1){1.1}}
\put(2.5,4){\line(-6,1){1.1}}
\put(3.3,2.9){\line(3,1){0.5}}
\put(1.7,2.9){\line(-3,1){0.5}}
\put(3.8,3.07){\line(-1,6){0.185}}
\put(1.2,3.07){\line(1,6){0.185}}
\put(3,2){\line(3,-4){0.34}}
\put(2,2){\line(-3,-4){0.34}}
\put(3.8,3.07){\line(1,-2){0.46}}
\put(1.2,3.07){\line(-1,-2){0.46}}
\put(3.35,1.55){\line(3,2){0.90}}
\put(1.65,1.55){\line(-3,2){0.90}}
\put(2.5,0.87){\line(-5,4){0.84}}
\put(2.5,0.87){\line(5,4){0.84}}
\put(2.5,0.87){\line(0,-1){0.69}}
\put(1,4.7){\line(1,0){3}}
\put(1,4.7){\line(3,-4){0.38}}
\put(4,4.7){\line(-3,-4){0.38}}
\put(1,4.7){\line(-1,-3){0.9}}
\put(4,4.7){\line(1,-3){0.9}}
\put(0.75,2.15){\line(-4,-1){0.65}}
\put(4.25,2.15){\line(4,-1){0.65}}
\put(2.5,0.18){\line(-4,3){2.4}}
\put(2.5,0.18){\line(4,3){2.4}}
\end{picture}
\caption{\textsf{do13} with five-fold rotational symmetry about the axis passing through the centers of faces $F_{11}$ and $F_{12}$}
\end{figure}

Let $c = \cos(\pi/5)$, $c_3 = \cos(3\pi/5)$, $s = \sin(\pi/5)$, $s_3 = \sin(3\pi/5)$, and
\begin{displaymath}
d = \sqrt{\frac{1+c}{1-c}} \frac{1}{\sqrt{2c+3c^2}+\sqrt{-1+3c^2}}.
\end{displaymath}
We obtain the following unit normals for the faces of \textsf{do13}:
{\setlength\arraycolsep{2pt}
\begin{eqnarray*}
\nu_1 &=& \big(\sqrt{ \frac{c}{4-4c}},\:  \frac{1}{\sqrt{2}s},\:0,\:\sqrt{\frac{2+3c}{4+4c}} \big),\\
\nu_{k+1} &=& L^{k}\nu_1 \:\:\text{with}\:\: k=1, \ldots, 9,\\
\nu_{11} &=& \big(\frac{1}{\sqrt{d^2-1}},\:0,\:0, \frac{d}{\sqrt{d^2-1}}\big),\\
\nu_{12} &=& \big(\frac{1}{\sqrt{d^2-1}},\:0,\:0, -\frac{d}{\sqrt{d^2-1}}\big).
\end{eqnarray*}}
where
\begin{displaymath}
L = \left( \begin{array}{cccc}
1 & 0 & 0 & 0\\
0 & c & -s & 0 \\
0 & s & c & 0 \\
0 & 0 & 0 & -1
\end{array} \right).
\end{displaymath}
To check these, we note that direct computation gives:
{\setlength\arraycolsep{2pt}
\begin{eqnarray*}
& & \langle \nu_i, \nu_i \rangle = 1  \:\:\text{for all}\:\: i=1, \ldots, 12,\\
& & \langle \nu_{11}, \nu_{2j-1} \rangle = \langle \nu_{12}, \nu_{2j} \rangle = -1/2 \:\:\text{for all}\:\: j=1, \ldots, 5,\\
& & \langle \nu_1, \nu_2 \rangle = \langle \nu_2, \nu_3 \rangle = \cdots = \langle \nu_{9}, \nu_{10} \rangle = \langle \nu_{10}, \nu_{1} \rangle = -1/2,\\
& & \langle \nu_1, \nu_3 \rangle = \langle \nu_3, \nu_5 \rangle = \cdots = \langle \nu_{7}, \nu_{9} \rangle = \langle \nu_{9}, \nu_{1} \rangle = 0,\\
& & \langle \nu_2, \nu_4 \rangle = \langle \nu_4, \nu_6 \rangle = \cdots = \langle \nu_{8}, \nu_{10} \rangle = \langle \nu_{10}, \nu_{2} \rangle = 0,
\end{eqnarray*}}
where $\langle \cdot\,,\cdot \rangle$ denotes the Lorentzian inner product.
Further, the cone defined by the inequalities $\langle X, \nu_i \rangle \leq 0$ intersects $\mathbb{H}^3$, since $X=(1,0,0,0)$ satisfies these inequalities.

We now want to solve Vinberg's equations (\ref{eqn:normalize})--(\ref{eqn:Vinberg2}).
We were unable to compute a Gr\"{o}bner basis for the corresponding ideal
using Mathematica since the coordinates of the linear functionals $\alpha_i=J\nu_i$ are complicated and the number of variables is very large.

Instead we look for solutions of the special form
{\setlength\arraycolsep{2pt}
\begin{eqnarray*}
b_1 &=& \big(b_{1,1},\: \frac{\sqrt{2}}{s},\:0,\:b_{1,4} \big) \:\:\text{with}\:\: \alpha_{1}(b_{1})=2,\\
b_2 &=& \big(b_{2,1},\: \frac{\sqrt{2}c}{s}, \sqrt{2},\:b_{2,4} \big) \:\:\text{with}\:\: \alpha_{2}(b_{2})=2,\\
b_{2k+1} &=& L^{2k}b_1 \:\:\text{and}\:\: b_{2k+2} = L^{2k}b_2 \:\:\text{with}\:\: k=1, \ldots, 4,\\
b_{11} &=& (b_{11,1},{0},{0},b_{11,4}) \:\:\text{with}\:\: \alpha_{11}(b_{11})=2,\\
b_{12} &=& (b_{12,1},{0},{0},b_{12,4}) \:\:\text{with}\:\: \alpha_{12}(b_{12})=2.
\end{eqnarray*}}
Here
we assume that $b_{11,2}=b_{11,3}=b_{12,2}=b_{12,3}=0$
to maintain the rotational symmetry of \textsf{do13}.
Therefore, we can choose $\{ b_{1,4}, b_{2,4}, b_{11,4}, b_{12,4}\}$ as  \emph{four} free variables.

Recall that $\alpha_i= J \nu_i$ for $i=1, \ldots, 12,$ and
$a_{ij}= \alpha_i(b_j)$. From direct calculations, we obtain
$$ a_{2k+1,2k+1}=a_{1,1}=2 \hspace{3mm} \text{and} \hspace{3mm} a_{2k+2,2k+2}=a_{2,2}=2,$$
for all $k=1, \ldots, 4$. Furthermore,
{\setlength\arraycolsep{2pt}
\begin{eqnarray*}
a_{2k-1,2k+1}&=&a_{1,3}=a_{1,1}+(a_{1,3}-a_{1,1})= 2+ \alpha_1(b_3-b_1)\\
&=& 2 + \frac{-1+\cos(2\pi/5)}{s^2}=0,\\
a_{2k+1,2k-1}&=&a_{3,1}=a_{1,1}+(a_{3,1}-a_{1,1})=0, \\
a_{2k,2k+2}&=&a_{2,4}=a_{2,2}+(a_{2,4}-a_{2,2})= 2+ \alpha_2(b_4-b_2)\\
{} &=& 2 + \frac{c}{\sqrt{2}s}\frac{\sqrt{2}(c_3-c)}{s} + \frac{s}{\sqrt{2}s}\frac{\sqrt{2}(s_3-s)}{s} = 0, \\
a_{2k+2,2k}&=&a_{4,2}=a_{2,2}+(a_{4,2}-a_{2,2})= 0,
\end{eqnarray*}}
for $k=0, 1, \ldots, 4$, where the indices are taken modulo $10$.

If the \emph{three} equations
\begin{displaymath}
a_{11,1} a_{1,11}=1,\hspace{2mm} a_{12,2} a_{2,12}=1 , \hspace{2mm} \text{and} \hspace{2mm} a_{1,2}a_{2,1}=1
\end{displaymath}
are satisfied, then so are all the Vinberg equations (\ref{eqn:normalize})--(\ref{eqn:Vinberg2}) since
{\setlength\arraycolsep{2pt}
\begin{eqnarray*}
a_{11,2k-1}=a_{11,1} \hspace{1mm} &\text{and}& \hspace{1mm} a_{2k-1,11}=a_{1,11},\\
a_{12,2k}=a_{12,2} \hspace{1mm} &\text{and}& \hspace{1mm} a_{2k,12}=a_{2,12},\\
a_{2k+1,2k}=a_{2k+1,2k+2}=a_{1,10}=a_{1,2} \hspace{1mm} &\text{and}& \hspace{1mm}
a_{2k,2k+1}=a_{2k+2,2k+1}=a_{10,1}=a_{2,1},
\end{eqnarray*}}
for $k=0, 1,\ldots, 4$ and indices modulo $10$.
It follows that the dimension of the local restricted deformation space is at least $1$.

On the other hand, we find that the dimension of the infinitesimal restricted deformation space of \textsf{do13} is exactly $1$ by exact computations using Mathematica.
This implies that the dimension of the local restricted deformation space is exactly equal to $1$.
(See the file do13.nb at \cite{Lee} for the detailed calculations.)
\end{example}

\begin{remark}\label{rem:geometric}
This work is similar to our work on prisms.
Both examples can be understood by letting the reflection points of the top and the bottom faces lie on the axis of rotational symmetry.
By choosing the reflection point arbitrarily for the top face, we see that the other faces are orderable
up to rotational symmetry. One can also take a quotient orbifold under the rotational symmetry and
obtain an ``orderable orbifold". Then the other reflection points can be chosen using this ordering.
Here, the geometry can be used to show explicitly how the reflection vectors change.

\end{remark}

Finally, we describe our method for studying the orbifolds \textsf{do1}-\textsf{do12}; this makes use of
numerical computations.

\begin{example} \label{example:do1}
We illustrate this for the orbifold  \textsf{do1}, shown in Figure 9.

\begin{figure}[h]
\centering
\setlength{\unitlength}{1cm}
\begin{picture}(5,5.5)
\put(2.3,2.55){\large $\boldsymbol{F_1}$}
\put(1.55,3.45){\large $\boldsymbol{F_2}$}
\put(1.15,2.25){\large $\boldsymbol{F_3}$}
\put(2.3,1.35){\large $\boldsymbol{F_4}$}
\put(3.40,2.25){\large $\boldsymbol{F_5}$}
\put(2.95,3.45){\large $\boldsymbol{F_6}$}
\put(0.65,3.1){\large $\boldsymbol{F_7}$}
\put(1.25,1.25){\large $\boldsymbol{F_8}$}
\put(3.35,1.25){\large $\boldsymbol{F_9}$}
\put(3.85,3.1){\large $\boldsymbol{F_{10}}$}
\put(2.25,4.25){\large $\boldsymbol{F_{11}}$}
\put(4.95,2.55){\large $\boldsymbol{F_{12}}$}
\put(2.75,2.95){$2$} 
\put(2.95,2.35){$3$}
\put(2.45,2.05){$2$}
\put(1.9,2.35){$3$}
\put(2.1,2.95){$3$}
\put(2.55,3.65){$2$} 
\put(3.4,3.0){$3$}
\put(3.2,1.75){$2$}
\put(1.65,1.75){$3$}
\put(1.45,3.0){$2$}
\put(3.0,4.10){$3$} 
\put(4.05,2.6){$2$}
\put(2.95,0.95){$3$}
\put(1.0,1.6){$2$}
\put(1.1,3.55){$3$}
\put(3.75,3.55){$3$} 
\put(3.8,1.6){$3$}
\put(1.9,0.95){$3$}
\put(0.8,2.6){$3$}
\put(1.85,4.10){$3$}
\put(3.6,4.35){$2$} 
\put(4.45,2.1){$3$}
\put(2.55,0.5){$2$}
\put(0.4,2.1){$3$}
\put(1.25,4.35){$2$}
\put(2.4,4.8){$3$} 
\put(4.5,3.3){$3$}
\put(3.8,0.9){$2$}
\put(1.0,0.9){$3$}
\put(0.35,3.3){$2$}
\put(2,2){\line(1,0){1}}
\put(2,2){\line(-1,3){0.3}}
\put(3,2){\line(1,3){0.3}}
\put(2.5,3.5){\line(4,-3){0.8}}
\put(2.5,3.5){\line(-4,-3){0.8}}
\put(2.5,3.5){\line(0,1){0.5}}
\put(2.5,4){\line(6,1){1.1}}
\put(2.5,4){\line(-6,1){1.1}}
\put(3.3,2.9){\line(3,1){0.5}}
\put(1.7,2.9){\line(-3,1){0.5}}
\put(3.8,3.07){\line(-1,6){0.185}}
\put(1.2,3.07){\line(1,6){0.185}}
\put(3,2){\line(3,-4){0.34}}
\put(2,2){\line(-3,-4){0.34}}
\put(3.8,3.07){\line(1,-2){0.46}}
\put(1.2,3.07){\line(-1,-2){0.46}}
\put(3.35,1.55){\line(3,2){0.90}}
\put(1.65,1.55){\line(-3,2){0.90}}
\put(2.5,0.87){\line(-5,4){0.84}}
\put(2.5,0.87){\line(5,4){0.84}}
\put(2.5,0.87){\line(0,-1){0.69}}
\put(1,4.7){\line(1,0){3}}
\put(1,4.7){\line(3,-4){0.38}}
\put(4,4.7){\line(-3,-4){0.38}}
\put(1,4.7){\line(-1,-3){0.9}}
\put(4,4.7){\line(1,-3){0.9}}
\put(0.75,2.15){\line(-4,-1){0.65}}
\put(4.25,2.15){\line(4,-1){0.65}}
\put(2.5,0.18){\line(-4,3){2.4}}
\put(2.5,0.18){\line(4,3){2.4}}
\end{picture}
\caption{\textsf{do1}}
\end{figure}

To find the unit normals $\nu_i$ for  \textsf{do1} we make some use of the results for \textsf{do13}.
We choose one vertex $F_{ijk}=F_i \cap F_j \cap F_k$ of \textsf{do1} whose adjacent edges have orders $(2,3,3)$.
For example, we choose the vertex $F_{123}$ of \textsf{do1}.
Let $c = \cos(\pi/5)$ and $s = \sin(\pi/5)$. We set
\begin{displaymath}
\nu_1 = \big(\frac{1}{\sqrt{d^2-1}},\:0,\:0, \frac{d}{\sqrt{d^2-1}}\big)\;\;\text{with}\;\; d = \sqrt{\frac{1+c}{1-c}} \frac{1}{\sqrt{2c+3c^2}+\sqrt{-1+3c^2}}.
\end{displaymath}
\begin{displaymath}
\nu_2 = \big(\sqrt{ \frac{c}{4-4c}},\:  \frac{1}{\sqrt{2}s},\:0,\:\sqrt{\frac{2+3c}{4+4c}} \big)\;\;\text{and}\;\; \nu_3=L^2 \nu_2,
\end{displaymath}
where \begin{displaymath}
L = \left( \begin{array}{cccc}
1 & 0 & 0 & 0\\
0 & c & -s & 0 \\
0 & s & c & 0 \\
0 & 0 & 0 & -1
\end{array} \right).
\end{displaymath}
From the calculations for \textsf{do13} we know that
\begin{displaymath}
\langle \nu_1, \nu_2 \rangle = \langle \nu_1, \nu_3 \rangle = -\frac{1}{2}\;\; \text{and}\;\; \langle \nu_2, \nu_3 \rangle=0.
\end{displaymath}
We first fix these normals $\nu_1$, $\nu_2$ and $\nu_3$,
 and solve the system of hyperbolic equations (\ref{eqn:hyp}) for
\textsf{do1} using the unit normals for \textsf{do13} as initial values for a numerical calculation by Mathematica.
This gives \emph{numerical} solutions for the unit normals of \textsf{do1}.

We then compute the rank of Jacobian $D$ of \textsf{do1} by numerical computations in Mathematica, maintaining 150 digits of precision in internal computations. This shows the rank  is $48$, and thus $D$ is of full rank.
\end{example}

To confirm these numerical computations,
we compute the singular value decomposition of $D$, and list the minimum
singular value in table 2.
(Note that the singular value decomposition
behaves very accurately in numerical matrix computations.)

For all $i=1,\ldots, 12$, we find that the Jacobian matrix $D$ of \textsf{do\,$i$} has rank $48$, and thus each \textsf{do\,$i$} is projectively rigid relative to the mirrors.
The numerical computations using Mathematica are also available from the webpage \cite{Lee}.


\begin{thebibliography}{30}
\addcontentsline{toc}{section}{References}

  \bibitem{Andreev1970}
    E.M. Andreev,
    \newblock {\em On convex polytopes in Lobachevskii spaces},
    \newblock Mat. Sbornik \textbf{81} (1970), 445--478.

  \bibitem{Andreev19702}
    E.M. Andreev,
    \newblock {\em On convex polytopes of finite volume in Lobachevskii space},
    \newblock Mat. Sbornik. \textbf{83} (1970), 256--260.

  \bibitem{Benoist2006}
    Y. Benoist,
    \newblock {\em Convexes divisibles IV},
    \newblock Invent. Math. \textbf{164} (2006), 249--278.

  \bibitem{BenoistNote}
    Y. Benoist,
    \newblock {\em Five lectures on lattices in semisimple Lie groups},
    \newblock Available at {http://www.math.u-psud.fr/\~{}benoist/}

  \bibitem{BHaf}
    M.R. Bridson and A. Haefliger,
    \newblock {\em Metric spaces of non-positive curvature},
    \newblock Springer, Berlin, 1999.

  \bibitem{CEG}
    R.D. Canary, D.B.A. Epstein and P. Green,
    \newblock {\em Notes on notes of Thurston},
    \newblock In: Epstein, D.B.A.(ed.) Analytical and geometric aspects of hyperbolic space, pp. 3--92, Cambridge Univ. Press, Cambridge, 1987.

    \bibitem{cdcr}
    S. Choi,
    \newblock {\em Convex decompositions of real projective surfaces I, II},
    \newblock J. Differential Geom. \textbf{40} (1994), 165--208, 239--283.

  \bibitem{Choi2004}
    S. Choi,
    \newblock {\em Geometric Structures on orbifolds and holonomy representations},
    \newblock Geom. Dedicata \textbf{104} (2004), 161--199.

  \bibitem{Choi2005}
    S. Choi and W.M. Goldman,
    \newblock {\em The deformation spaces of convex $\mathbb{RP}^2$-structures on 2-orbifolds},
    \newblock Amer. J. Math. \textbf{127} (2005), 1019--1102.

  \bibitem{Choi2006}
    S. Choi,
    \newblock {\em The deformation spaces of projective structures on 3-dimensional Coxeter orbifolds},
    \newblock Geom. Dedicata \textbf{119} (2006), 69--90.

   \bibitem{Choi2010}
    S. Choi,
    \newblock {\em Convexity while deforming convex real projective manifolds and orbifolds with ends},
    \newblock Preprint (2010)

  \bibitem{Cooper2006}
    D. Cooper, D. Long and M. Thistlethwaite,
    \newblock {\em Computing varieties of representations of hyperbolic 3-manifolds into ${\rm SL}(4,\mathbb R)$},
    \newblock Experiment. Math. \textbf{15} (2006), 291--305.

 \bibitem{GarlandRag}
    H. Garland and M.S. Raghunathan,
    \newblock {\em Fundamental domains for lattices in (R-)rank $1$ semisimple Lie groups},
    \newblock Ann. of Math. \textbf{92} (1970), 279--326.

  \bibitem{Goldman1990}
    W.M. Goldman,
    \newblock {\em Convex real projective structures on compact surfaces},
    \newblock J. Differential Geom. \textbf{31} (1990), 791--845.

  \bibitem{orb} D. Heard,
    \newblock ``Orb'':  a computer program for finding
    hyperbolic structures on hyperbolic 3-orbifolds and $3$-manifolds.
    \newblock Available at {http://www.ms.unimelb.edu.au/\~{}snap/orb.html}

  \bibitem{Lee}
    G.-S. Lee,
    \newblock Matlab and Mathematica files for the computations in this paper.
    \newblock Available at {http://mathsci.kaist.ac.kr/\~{}manifold}

 \bibitem{Kapovich}
     M. Kapovich,
    \newblock {\em Hyperbolic manifolds and discrete groups},
     \newblock Birkh\a"user, Boston, 2001.

 \bibitem{Kapovich1994}
    M. Kapovich,
    \newblock {\em Deformations of representations of discrete subgroups of ${\rm SO}(3,1)$},
    \newblock Math. Ann. \textbf{299} (1994), 341--354.

    \bibitem{Raghunathan}
     M.S. Raghunathan,
     \newblock {\em Discrete subgroups of Lie groups},
     \newblock Springer, Berlin, 1972.

  \bibitem{Roeder2007}
    R.K.W. Roeder,
    \newblock {\em Constructing hyperbolic polyhedra using Newton's method},
    \newblock Experiment. Math. \textbf{16} (2007), 463--492.

  \bibitem{RHD}
    R.K.W. Roeder, J.H. Hubbard and W.D. Dunbar,
    \newblock {\em Andreev's theorem on hyperbolic polyhedra},
    \newblock Ann. Inst. Fourier (Grenoble) \textbf{57} (2007), 825--882.

  \bibitem{Thurston}
    W. Thurston,
    \newblock {\em Geometry and topology of 3-manifolds},
    \newblock Lecture notes. Princeton University (1979)
    \newblock Available at http://www.msri.org/publications/books/gt3m/

    \bibitem{Thurston_book}
    W. Thurston,
    \newblock {\em Three-dimensional geometry and topology},
    \newblock Princeton University Press, New Jersey, 1997.

  \bibitem{Trefethen}
    L.N. Trefethen and D. Bau, III:
    \newblock {\em Numerical Linear Algebra},
    \newblock SIAM, Philadelphia, 1997.

  \bibitem{Kac1967}
    \`{E}.B. Vinberg and V.G. Kac,
    \newblock {\em Quasi-homogeneous cones},
    \newblock Math. Zamnetki \textbf{1} (1967), 347--354.

  \bibitem{Vinberg1971}
    \`{E}.B. Vinberg,
    \newblock {\em Discrete linear groups that are generated by reflections},
    \newblock Izv. Akad. Nauk SSSR Ser. Mat. \textbf{35} (1971), 1072--1112.

  \bibitem{Vinberg1985}
    \`{E}.B. Vinberg,
    \newblock {\em Hyperbolic reflection groups},
    \newblock Uspekhi Mat. Nauk \textbf{40} (1985), 29--66.

  \bibitem{Vinberg1993}
    \`{E}.B. Vinberg(ed.),
    \newblock {\em Geometry II},
    \newblock Springer, Berlin, 1993.

  \bibitem{Weil1962}
    A. Weil,
    \newblock {\em On Discrete subgroups of Lie groups II},
    \newblock Ann. of Math. \textbf{75} (1962), 578--602.

   \bibitem{Weil1964}
    A. Weil,
    \newblock {\em Remarks on the cohomology of groups},
    \newblock Ann. of Math. \textbf{80} (1964), 149--157.


\end{thebibliography}
\end{document}